\documentclass[12pt]{article}

\usepackage{multirow}
\usepackage{color}

\usepackage[a4paper, margin=0.8in]{geometry}

\usepackage{booktabs}
\usepackage{tabularx}
\usepackage{amssymb}
\usepackage{gensymb}
\newcolumntype{Y}{>{\centering\arraybackslash}X}
\usepackage{setspace}
\usepackage{graphicx}
\usepackage[utf8]{inputenc}
\usepackage{subfigure}
\usepackage{rotating}
\usepackage{amsmath}
\usepackage{amssymb}
\usepackage{float}
\usepackage{bm}
\usepackage[binary-units]{siunitx}
\usepackage[left,modulo]{lineno}
\usepackage[colorlinks=true]{hyperref}
\usepackage[affil-it]{authblk}
\usepackage[square,sort,numbers]{natbib}

\newcommand{\bn}{{\bm n}}
\newcommand{\bx}{{\bm x}}
\newcommand{\by}{{\bm y}}

\newcommand{\rmi}{{\rm i}}

\providecommand{\keywords}[1]{\textbf{\textit{Keywords---}} #1}

\title{A-priori mesh grading for the numerical calculation\\of the head-related transfer functions\footnote{Manuscript submitted to Applied Acoustics.}}
\author[1]{Harald Ziegelwanger}
\author[2]{Wolfgang Kreuzer}
\author[2]{Piotr Majdak}

\affil[1]{AIT Austrian Institute of Technology GmbH, Mobility Department, Transportation Infrastructure Technologies, Giefinggasse 2, 1210 Vienna, Austria.\\ email: harald.ziegelwanger@ait.ac.at}
\affil[2]{Acoustics Research Institute, Austrian Academy of Sciences, Wohllebengasse 12-14, Vienna, 1040, Austria}

\date{}

\begin{document}
\maketitle

\abstract{
    Head-related transfer functions (HRTFs) describe the directional filtering of the incoming sound caused by the morphology of a listener's head and pinnae. When an accurate model of a listener's morphology exists, HRTFs can be calculated numerically with the boundary element method (BEM). However, the general recommendation to model the head and pinnae with at least six elements per wavelength renders the BEM as a time-consuming procedure when calculating HRTFs for the full audible frequency range.
    In this study, a mesh preprocessing algorithm is proposed, viz., a-priori mesh grading, which reduces the computational costs in the HRTF calculation process significantly.
    The mesh grading algorithm deliberately violates the recommendation of at least six elements per wavelength in certain regions of the head and pinnae and varies the size of elements gradually according to an a-priori defined grading function.
    The evaluation of the algorithm involved HRTFs calculated for various geometric objects including meshes of three human listeners and various grading functions.
    The numerical accuracy and the predicted sound-localization performance of calculated HRTFs were analyzed.
    A-priori mesh grading appeared to be suitable for the numerical calculation of HRTFs in the full audible frequency range and outperformed uniform meshes in terms of numerical errors, perception based predictions of sound-localization performance, and computational costs.
}




\section{Introduction}
\label{sec:Introduction}
The head-related transfer functions (HRTFs) describe the directional filtering of incident sound waves at the entrance of a listener's ear-canal \citep{mehrgardt_transformation_1977,wightman_headphone_1989}. This filtering is caused by the overall diffraction, shadowing, and reflections at human anatomical structures, i.e., the torso, head, and pinnae. Thus, HRTFs are closely related to a listener's individual geometry and they provide listener-specific spectral \citep{middlebrooks_individual_1999} and temporal features \citep{macpherson_listener_2002} which are essential for three-dimensional (3D) sound localization, e.g., in binaural audio reproduction systems \citep{moller_fundamentals_1992}.

HRTFs are usually acquired acoustically in a resource-demanding process, in which small microphones are placed into listener's ear canals and transfer functions are measured for many directions in an anechoic chamber \citep{moller_head-related_1995,majdak_3-d_2010}. HRTFs can also be acquired by means of a numerical HRTF calculation, i.e., by simulating the sound field of an incident wave scattered by a listener's head and pinnae. In recent years, the boundary element method \citep[BEM,][]{gaul_boundary_2003} became a powerful tool for such simulations in acoustics. The BEM was used in many studies for the numerical calculation of HRTFs \citep{kahana_numerical_2006,kahana_boundary_2007,gumerov_computation_2010,katz_boundary_2001,katz_boundary_2001-1,kreuzer_fast_2009,rui_calculation_2013,jin_creating_2014,ziegelwanger_numerical_2015}. In general, the numerical HRTF calculation is based on a 3D polygon mesh, i.e., a set of nodes and elements with these nodes as vertices, describing a listener's morphology.

In element-based acoustic simulations resolution of the mesh should be related to the wavelength of the simulated frequency \citep{marburg_six_2002}. The mesh resolution is measured by the number of elements per wavelength or by the average length of edges in the mesh, i.e., the average edge length \citep[AEL,][]{ziegelwanger_numerical_2015}. The accuracy of the numerical calculations depends on these metrics. In \citet{marburg_six_2002} the relative numerical error was below fifteen percent, when at least six elements per wavelength were used. \citet{gumerov_computation_2010} recommended five elements per wavelength, equilateral triangles, and a valence of six, i.e., the number of edges incident to a vertex describing the regularity of a mesh, and a uniform vertex distribution in the mesh. In \citet{ziegelwanger_numerical_2015}, an AEL of \SIrange{1}{2}{\milli\meter} was required for accurate numerical HRTF calculations. Given the average human body surface area and a frequency range of up to \SI{18}{\kilo\hertz}, these recommendations result in a uniform head and pinna mesh consisting of approximately \num{100000} equilateral triangular elements.

The computational costs of the BEM, i.e., the calculation time and the required amount of physical memory, are generally high and increase with the number of elements in the mesh \citep{kreuzer_fast_2009}. The first numerical HRTF calculations were limited to \num{22000} elements and frequencies up to \SI{5,4}{\kilo\hertz} because the calculation time was in the range of tens of hours for a single frequency \citep{katz_boundary_2001}. The HRTF calculation became feasible for the full audible frequency range \citep{kreuzer_fast_2009,gumerov_computation_2010} by coupling the BEM with the fast multipole method \citep[FMM, e.g.,][]{chen_formulation_2008}. HRTFs calculated with the FMM showed good results for artificial heads by means of visual comparison of amplitude spectra \citep{gumerov_computation_2010} and for human listeners by means of individual sound-localization performance \citep{ziegelwanger_numerical_2015}. However, the numerical HRTF calculation process for the full audible frequency range can still last tens of hours on a single CPU \citep{kreuzer_fast_2009,ziegelwanger_numerical_2015}.

While the computational costs can be reduced by reducing the number of elements in the mesh, a simple mesh coarsening considering all elements, i.e., a \emph{uniform} re-meshing, also reduces the accuracy of the numerical calculation \citep{marburg_six_2002}. The loss of accuracy is caused by geometric and numerical error \citep{treeby_practical_2009}. The geometric error arises because of under-sampling the geometry and the numerical error arises because of under-sampling the sound field on the geometry \citep{ziegelwanger_numerical_2015}. In other fields of computational physics, more sophisticated geometry discretization methods, resulting in \emph{non-uniform} meshes, have been investigated. For the finite-element method, the numerical error introduced by goal-oriented mesh adaptivity \citep{walsh_hp_2003}, non-uniform meshes \citep{goldstein_finite_1982}, and mesh grading \citep{heinrich_optimum_1996} were investigated. For elliptic boundary-value problems, \emph{a-priori mesh grading} was proposed \citep{Langer20151685}, where the element size was varied based on a-priori knowledge of the location of singular points, i.e., sharp edges and corners in the geometry or discontinuities in the boundary conditions. For the 2D-BEM and hyperbolic boundary-value problems, adaptive meshes were investigated \citep{chen_adaptive_2002,liang_error_1999}. In general, non-uniform meshes showed better convergence rates than uniform meshes and the accuracy was higher for non-uniform meshes than for uniform meshes. Even though the investigations for hyperbolic boundary-value problems were done for the 2D-BEM only, a non-uniform mesh of a human head seems to be a promising approach to reduce the computational costs in the numerical HRTF calculation process.

Hence, in this article, we adapt the idea of a-priori mesh grading for the numerical calculation of HRTFs. First, we briefly review the BEM (Sec.~\ref{sec:Boundary-element-method}) and describe our a-priori mesh grading algorithm (Sec.~\ref{sec:Mesh-grading}). In Sec.~\ref{sec:Evaluation}, we show its evaluation. Secs.~\ref{sub:HRTFcalculation} and \ref{sub:Metrics} describe the software and metrics we have used for the evaluation. Sec.~\ref{sub:Experiment11} and \ref{sub:Experiment12}, show the evaluation of various grading functions based on a comparison to reference HRTFs of a sphere. Then, the most promising grading functions were evaluated on the geometry of a pinna (Sec.~\ref{sub:Experiment2}). Finally, in Sec.~\ref{sub:Experiment3}, the best performing grading function was applied on meshes of human heads, for which the HRTFs were evaluated by means of numerical and perceptual errors.


\section{Boundary element method}
\label{sec:Boundary-element-method}
\begin{figure}[!h]
    \centering
    \includegraphics[width=0.6\textwidth]{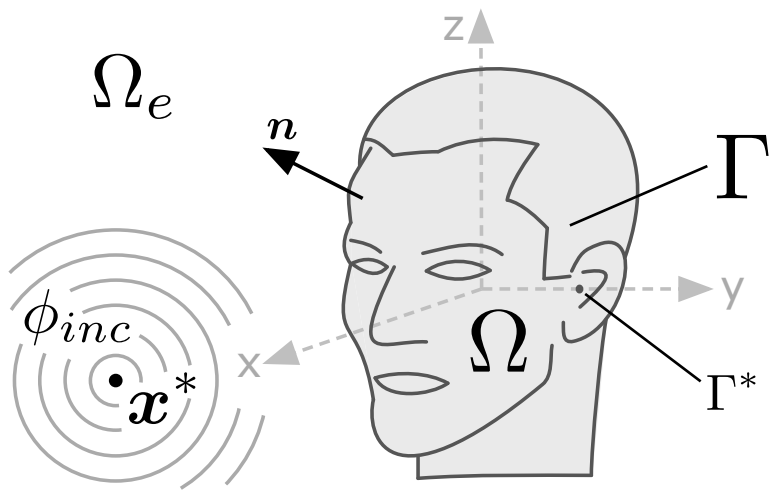}
    \caption{Schematic representation of the exterior scattering problem for the numerical calculation of HRTFs. A point source is placed at $\bx^*$ and emits the incident sound field $\phi_{inc}(\bx)$ in $\Omega_e$ outside a listeners head $\Omega$ with surface $\Gamma$ and the microphone area $\Gamma^*$. x, y, and z represent the Cartesian coordinate system as described in \citet{ziegelwanger_modeling_2014}.}
\label{fig:BEM}
\end{figure}

The boundary element method in the context of calculating HRTFs is schematically shown in Fig.~\ref{fig:BEM}. Here, the object $\Omega$ with boundary $\Gamma$ represents the scatterer, i.e. the human head and pinnae. $\Omega_e$ defines the domain outside the scatterer. A point source at $\bx^*$ (in the following called \emph{loudspeaker position}) emits spherical waves, i.e., it produces the incident sound field $\phi_{inc}(\bx)$. $\Gamma^*$ is the \emph{microphone area} at the entrance of the ear canal.

The total sound field $\phi(\bx)$ is described by the boundary integral equation:
\begin{equation}
  \tau \phi(\bx)\!=\!\int\limits_{\Gamma} H(\bx,\by) \phi(\by)d\by \,  - \int\limits_{\Gamma} G(\bx,\by) v(\by) d\by \, + \phi_{inc}(\bx),
  \label{equ:BIE}
\end{equation}
where $G(\bx,\by)$ and $H(\bx,\by)$ are the Green's function of the Helmholtz equation and its derivative with respect to the normal vector $\bn$ to $\Gamma$ at the point $\by$.
$\phi(\bx) = \frac{p(\bx)}{\rmi \omega \rho}$, $p(\bx)$ and $v(\bx) = \frac{\partial \phi(\bx)}{\partial {\bn}} = {\bn} \cdot \nabla \phi(\bx)$ denote the velocity potential, the sound pressure, and the particle velocity at a point $\bx$, respectively.
$\rho$ denotes the density of air and $\tau$ is a factor depending on the position of $\bx$. $\tau = 1$ for $\bx \in \Omega_e$, $\tau = 1/2$ for $\bx \in \Gamma$, and $\tau=0$ when $\bx$ is located inside $\Omega$. The scatterer is assumed to be rigid, thus $\frac{\partial \phi(\bx)}{\partial \bn} = 0$ for $\bx \in \Gamma$.

To speed up calculations, HRTFs are determined by applying the principle of reciprocity \citep{morse_theoretical_1986}, where the roles of sources and receivers are exchanged. To this end, the in-ear microphone is simulated by a point source close to $\Gamma^*$ \citep{gumerov_computation_2010} or by active vibrating elements at $\Gamma^*$ \citep{ziegelwanger_numerical_2015,kreuzer_fast_2009}. HRTFs are evaluated using the calculated sound pressure at the loudspeaker positions. For the active vibrating elements, this means technically that the contribution of an external sound source $\phi_{inc}(\bx)$ is replaced by an additional boundary condition $v(\bx) \neq 0$ for $\bx \in \Gamma^*$.

In our approach, HRTFs are calculated numerically in three steps. First, $\Gamma$ is discretized as a 3D polygon mesh $\mathcal{M}$, consisting of vertices $\mathcal{V}$, edges $\mathcal{E}$ and elements $\mathcal{F}$ (Fig.~\ref{fig:meshOperationsRemeshing}a), and the unknown solution $\phi(\bx)$ on $\Gamma$ is approximated using simple basis functions \citep{hunter_fem/bem_2002}, e.g., piecewise constant basis functions.
Using a collocation approach Eq.~\ref{equ:BIE} (for $\tau=0.5$) is transformed into a linear system of equations $\mathbf{A} \boldsymbol{\phi}=\mathbf{b}$, where $\boldsymbol{\phi}$ is the vector of unknown velocity potentials. The Burton-Miller approach \citep{burton_application_1971} is used to ensure a unique solution. For details about the derivation of the stiffness matrix $\mathbf{A}$ and the right-hand-side $\mathbf{b}$ refer to \citet{chen_formulation_2008} and \citet{ziegelwanger_mesh2hrtf:_2015}.
Second, the solution for the linear system of equations is calculated by using an iterative solver. The FMM is used to speed up the matrix-vector multiplications needed for the iterative solver \citep{chen_formulation_2008}. Third, given the solution $\boldsymbol{\phi}$ at the boundary $\Gamma$, the sound pressure $p(\bx)=\rmi \rho \omega \phi(\bx)$ at any point $\bx$ in the exterior domain $\Omega_e$, e.g., the loudspeaker grid, is determined by evaluating Eq.~\ref{equ:BIE} (for $\tau=1$).


\section{A-priori mesh grading}
\label{sec:Mesh-grading}
The a-priori mesh grading approach aims at reducing the number of elements $\#\mathcal{F}$ while preserving the accuracy in the calculation results by gradually increasing the length of edges in $\mathcal{M}$ as a function of the distance of an edge to $\Gamma^*$. The validity of our approach is based on two assumptions.

As for the geometric error, we assume that the geometry of the ipsilateral pinna\footnote{When using the reciprocity approach, HRTFs need to be calculated for both ears separately. Thus, in this context, the terms ipsilateral and contralateral correspond to the side of the head, where the sound source is located or not, respectively.} has the main influence on the accuracy of the HRTFs \citep[see also][]{ziegelwanger_numerical_2015}, whereas the geometry of the rest of the head (including the contralateral pinna) plays a minor role. Thus, the effect of a geometric error caused by large elements at the contralateral head side is assumed to be negligible.

As for the numerical error, we assume that the discretization of $\phi(\bx)$ in the proximity range of a discontinuity in boundary conditions, i.e., the jump in the neumann boundary condition at $\Gamma^*$ in the reciprocal calculation, has to be very fine, whereas the effect of the size of elements in regions far apart from $\Gamma^*$ is negligible, because the Green's function $G(\bx,\by)$ and its derivatives in Eq.~\ref{equ:BIE} decay fast with increasing distance $||\bx-\by||$. Additionally the amplitude of $\phi(\bx)$ on the contralateral side will be much smaller than on the ipsilateral side (with a difference of approximately \SI{90}{\dB}). Therefore, any error caused by the increased element size at the contralateral side will be damped in the order of tens of \si{\dB}, introducing only small artifacts to the final HRTFs.

Based on our assumptions, an algorithm which gradually increases the length of edges in $\mathcal{M}$ as a function of the distance of an edge to $\Gamma^*$ was designed. The algorithm is based on a re-meshing algorithm for uniform meshes \citep{botsch_remeshing_2004}. In this algorithm, if the length of an edge does not match a pre-defined target edge length, the mesh around this edge is modified. In our approach, this algorithm is extended such that the target edge length is described by a grading function dependent on the (relative) distance of an edge to $\Gamma^*$. In the following, we describe the details of the grading functions and of the re-meshing algorithm.
\subsection{Grading functions}
\label{sub:Grading-functions}
We consider the relative distance of an edge $\mathcal{E}$ to $\Gamma^*$:
\begin{equation}
        \bar{d}_{\mathcal{E}}=\frac{| \mathcal{E}_m-\Gamma_m^* |}{d_{max}},
        \label{equ:NormalizedDistance}
\end{equation}
where $\mathcal{E}_m$ is the midpoint of $\mathcal{E}$, $\Gamma_m^*$ is the midpoint of $\Gamma^*$, and $d_{max} = {\max_{\mathcal{E}}}\, d_{\mathcal{E}}$ is the maximum distance of edges in $\mathcal{M}$.
Note that $\bar{d}_{\mathcal{E}}=0$ and $\bar{d}_{\mathcal{E}}=1$ correspond to the most ipsilateral and most contralateral side of the head, respectively.

\begin{figure}[!h]
    \centering
    \includegraphics[width=0.6\linewidth]{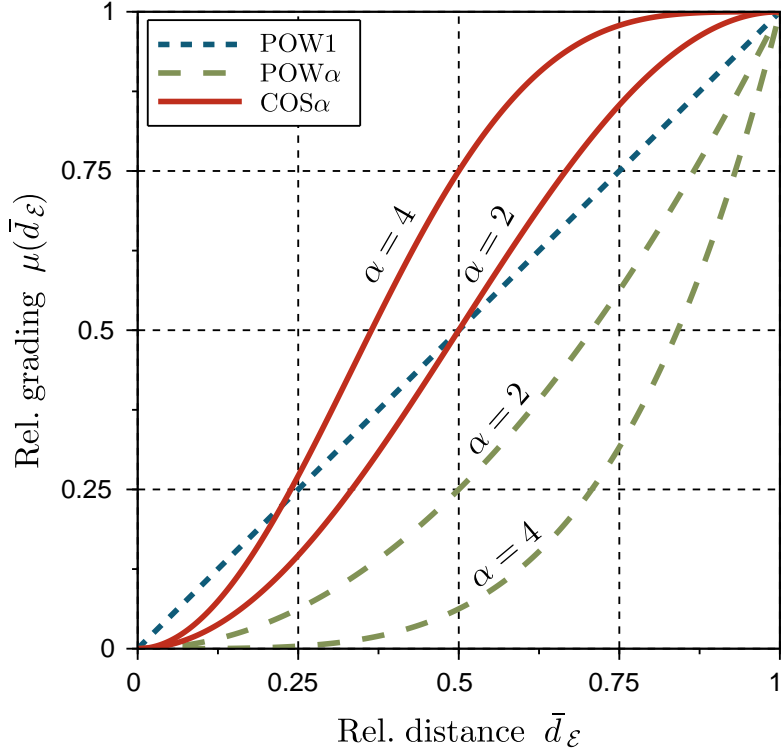}
    \caption{Relative grading resulting from various grading functions (linear POW1, higher-order grading POW$\alpha$, and raised-cosine grading COS$\alpha$).}
\label{fig:gradingFunctions}
\end{figure}

For the grading of $\mathcal{M}$, we propose various grading functions $\mu(\bar{d}_{\mathcal{E}})$. These functions, shown in Fig.~\ref{fig:gradingFunctions}, can be structured in two classes.
The first class consists of power functions:
\begin{equation}
	\mu(\bar{d}_{\mathcal{E}})={\bar{d}_{\mathcal{E}}}^{\alpha},
\label{equ:QuadraticGrading}
\end{equation}
with $\alpha=1$ as the first-order grading (\textbf{POW1}), which increases the length of edges linearly from the ipsilateral to the contralateral head side. Also, higher-order grading functions were considered, e.g., quadratic grading (\textbf{POW2}) and quartic grading (\textbf{POW4}) for $\alpha=2$ and $\alpha=4$, respectively. With increasing $\alpha$, these grading functions increase the element size slightly at the ipsilateral side and rapidly at the contralateral side. Note that a zeroth-order grading yields a uniform mesh.

The second class is based on raised cosine:
\begin{equation}
  \mu(\bar{d}_{\mathcal{E}})=1-\cos^{\alpha}\left(\pi \bar{d}_{\mathcal{E}}/2\right),
\label{equ:CosineGrading}
\end{equation}
\noindent with the second-order cosine grading (\textbf{COS2}) for $\alpha=2$ and fourth-order cosine grading (\textbf{COS4}) for $\alpha=4$. In contrast to the POW$\alpha$ functions, COS$\alpha$ functions concentrate small edges around $\Gamma^*$. While $\alpha=2$ then increases the edge length almost linearly from the ipsilateral to the contralateral head side, $\alpha=4$ grading increases the edge length rapidly to the maximum edge length.

In order to calculate the target edge length for each edge in $\mathcal{M}$, two global mesh grading parameters are defined, the minimum target edge length $\hat{\ell}_{min}$ and the maximum target edge length $\hat{\ell}_{max}$. The local target edge length $\hat{\ell}_{\bar{d}_\mathcal{E}}$ is then calculated as:

\begin{equation}
    \hat{\ell}_{\bar{d}_\mathcal{E}} = \hat{\ell}_{min} + \left ( \hat{\ell}_{max}-\hat{\ell}_{min} \right ) \mu(\bar{d}_{\mathcal{E}}).
    \label{equ:TargetEdgeLength}
\end{equation}

\subsection{Re-meshing algorithm}
\label{sub:Re-meshing}
\begin{figure}[!h]
    \centering
    \includegraphics[width=0.5\linewidth]{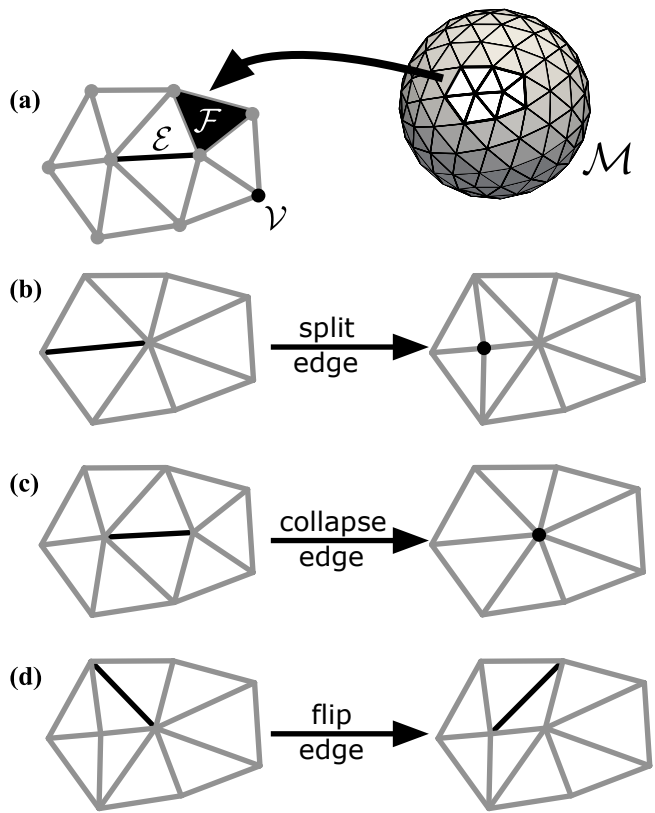}
    \caption{The re-meshing algorithm. \textbf{(a)} Definitions: $\mathcal{V}$, $\mathcal{E}$, and $\mathcal{F}$ represent a vertex, an edge, and a face, respectively. Mesh operations: \textbf{(b)} edge splitting, \textbf{(c)} edge collapsing, and \textbf{(d)} edge flipping.}
\label{fig:meshOperationsRemeshing}
\end{figure}

Most of the re-meshing algorithms in the field of computer graphics use triangle mesh decimation, vertex clustering, or voxel based object simplification \citep{luebke_developers_2001}. Our a-priori mesh grading algorithm is based on an efficient re-meshing algorithm from \citet{botsch_remeshing_2004}\footnote{Available in OpenFlipper, version 2.1,  \citep{mobius_openflipper:_2012}: \url{http://www.openflipper.org} (date last viewed: January 31, 2016)}. This algorithm modifies a mesh in an iterative procedure with the goal to obtain a uniform mesh with a given target edge length. In our mesh grading algorithm, the target edge length is the local target edge length $\hat{\ell}_{\bar{d}_\mathcal{E}}$ from the previous section.
Thus, in each iteration, first, edges are split if $\ell_{\mathcal{E}} > \frac{4}{3} \hat{\ell}_{\bar{d}_\mathcal{E}}$ (Fig.~\ref{fig:meshOperationsRemeshing}b).
Second, edges are collapsed if $\ell_{\mathcal{E}} < \frac{4}{5} \hat{\ell}_{\bar{d}_\mathcal{E}}$ (Fig.~\ref{fig:meshOperationsRemeshing}c).
Third, edges are flipped if the valence of neighboring vertices is larger than six (Fig.~\ref{fig:meshOperationsRemeshing}d).
Last, vertices are relocated on the surface of the original mesh by tangential smoothing \citep{botsch_remeshing_2004}. In our study, ten iterations were sufficient to achieve the edge length distribution targeted by the corresponding grading function.

\section{Evaluation}
\label{sec:Evaluation}
The effect of a-priori mesh grading on numerically calculated HRTFs was investigated in three steps, with different geometric objects in each step. 
First, HRTFs were numerically calculated for a rigid sphere (\textbf{SPH}) object and were compared to an analytical solution. Second, the same HRTFs were compared to HRTFs calculated for a uniform high-resolution mesh of that sphere, validating the usability of a high-resolution mesh to generate a reference HRTF set for further investigations. Third, the SPH object was extended by a generic pinna yielding the sphere-and-pinna (\textbf{SAP}) object. For that object, an analytical solution is not available and HRTFs numerically calculated from a high-resolution mesh were used as the reference. Fourth, HRTFs were calculated for three human objects (\textbf{HUMs}) represented by head and pinna models of actual human listeners. For the HUMs objects, the effect of a-priori mesh grading on numerically calculated HRTFs was (in addition to numerical errors) evaluated from the perceptual point of view. In this section, we first describe some general aspects of our evaluation and then show and discuss the results for the evaluated objects.

\subsection{HRTF calculation}
\label{sub:HRTFcalculation}
Mesh2HRTF\footnote{Mesh2HRTF: version 0.1.2, available from \url{http://mesh2hrtf.sourceforge.net} (date last viewed: January 31, 2016)} was used to calculate HRTFs numerically for \num{200} frequencies, which were linearly spaced between \SIlist{0,1;20}{\kilo\hertz}. Calculations were run on a Linux cluster consisting of eight PCs with Intel i7--3820 processors running at \SI{3.6}{\giga\hertz}. Each PC was equipped with \SI{64}{\giga\byte} of RAM. In total, more than hundred HRTF sets were calculated in this study and the calculations lasted several days.

HRTF positions were defined by means of two loudspeaker grids. A loudspeaker grid represents the evaluation nodes in the BEM. In the first grid, the equi-angular (\textbf{EQA}) grid, the interaural-polar coordinate system was used \citep{morimoto_localization_1984}. In that system, the poles are placed along the interaural axis and a direction is represented by lateral angle $\vartheta$ and polar angle $\varphi$. In that grid, \num{5042} virtual loudspeakers were placed equi-angular on a sphere with radius \SI{1,2}{\meter}. The resolution was \ang{2,5} and \ang{5} in the lateral and polar dimension, respectively. The EQA grid was used to analyze the \emph{numerical} accuracy of the numerically calculated HRTFs.

In the second grid, the \textbf{ARI} grid, \num{1550} virtual loudspeakers were placed on a sphere with radius \SI{1,2}{\meter} \citep{majdak_3-d_2010,ziegelwanger_numerical_2015}. The spherical coordinate system was used which is described by the azimuth and elevation angle. The grid had a polar gap at elevation angles below \ang{-30}, as well as a resolution of \ang{5} in the elevation and \ang{2,5} and \ang{5} for azimuth angles smaller and larger than \ang{30}, respectively. The ARI grid was used in the \emph{perceptual} evaluation of the numerically calculated HRTFs.

\subsection{Considered metrics}
\label{sub:Metrics}
Two error metrics were used in the evaluation. First, the accuracy of numerical calculation was considered which shows the performance and limits of the method and frames our results into context of previous studies in the field of computational acoustics. Second, individual sound-localization performance was considered which shows the perceptual impact of the numerical error and the relevance of the method in the context of HRTFs and frames our research into context of previous studies in the field of spatial hearing.

The accuracy of the numerical calculation was quantified by the relative numerical error. To this end, the complex-valued sound pressure from a numerical calculation was compared with a reference solution resulting in the relative numerical error $e^{\Omega_e}_{L^{\mathrm{p}}}$ given by:
\begin{equation}
    e^{\Omega_e}_{L^{\mathrm{p}}}=\frac{\lVert p_{num}^{\Omega_e} - p_{ref}^{\Omega_e} \rVert_{L^{\mathrm{p}}}}{\lVert p_{ref}^{\Omega_e} \rVert_{L^{\mathrm{p}}}},
    \label{equ:NumericalError}
\end{equation}
\noindent where $p_{num}$ and $p_{ref}$ are the evaluated and a reference sound pressure, respectively. $||.||_{L^{\mathrm{p}}}$ denotes an $L^{\mathrm{p}}$-norm.

Two $L^{\mathrm{p}}$-norms were considered, viz., the $L^{\mathrm{2}}$-norm and the $L^{\mathrm{\infty}}$-norm. For the calculated frequencies and the EQA grid, the $L^2$ norm was calculated as:
\begin{equation}
    || f ||_{L^2} \approx \sqrt{ \Delta \omega \Delta \varphi \Delta \vartheta \sum_i \sum_j \sum_k | f (\omega_i, \varphi_j, \vartheta_k)|^2 \sin(\vartheta_k) },
    \label{equ:NormSphere2}
\end{equation}
where $(\varphi_j, \vartheta_k)$ are defined by the EQA grid. The $L^{\mathrm{\infty}}$-norm was:
\begin{equation}
    \lVert f \rVert_{L^{\infty}}=\max_{\omega,\varphi,\vartheta} | f(\omega,\varphi,\vartheta) |.
    \label{equ:NormSphereInf}
\end{equation}
Note that our error formulation in Eq.~\ref{equ:NumericalError} does not separate magnitude and phase differences and thus considers both. A phase difference of $\pi$ results in a relative error of two hundred percent.

The individual sound-localization performance was quantified in terms of model predictions for an HRTF set. Two models were used. First, spectral features were analyzed with a model of sound-localization performance in sagittal-planes \citep{baumgartner_assessment_2013} implemented as \texttt{baumgartner2013} in the Auditory Modeling Toolbox\footnote{AMT: version 0.9.6, available from \url{http://amtoolbox.sourceforge.net} (date last viewed: January 31, 2016)} \citep[AMT,][]{sondergaardmajdak2013}. This model predicts the individual sound-localization performance by means of the local polar RMS error (PE) and the quadrant error rate (QE), which are common error metrics for analyzing results of sound-localization experiments \citet{middlebrooks_individual_1999}. Second, temporal features were analyzed by a direction-continuous model of the time-of-arrival \citep[TOA,][]{ziegelwanger_modeling_2014} implemented as \texttt{ziegelwanger2014} in the AMT. In the TOA model, the temporal features are quantified as the equivalent head radius resulting from the tested HRTF set. A deviation from the radius obtained for a reference HRTF set can be interpreted as an artifact in the broadband temporal properties of tested HRTFs, having potentially an effect on the perception of the interaural time differences with such HRTFs.

\subsection{The SPH object and analytical reference}
\label{sub:Experiment11}
The SPH object was a rigid sphere model with a radius $R$ of \SI{100}{\milli\meter}. First, a high-resolution mesh of the sphere was constructed in Blender\footnote{Blender: version 2.72b, available from \url{http://www.blender.org} (date last viewed: January 31, 2016)}. It consisted of \num{139194}~elements corresponding to a resolution of about \num{14} elements per wave length at \SI{18}{\kilo\hertz} and an AEL of approximately \SI{1,4}{\milli\meter}. Then, by uniformly re-meshing that mesh, uniform meshes (UNI) were created with AELs ranging from \SIrange{2}{20}{\milli\meter}. Further, by applying the a-priori mesh grading to the high-resolution mesh, graded meshes were created with $\#\mathcal{F}$s and AELs covering the range of the uniform meshes. Five grading functions were considered (POW1, POW2, POW4, COS2, and COS4). Table~\ref{tab:sphereMeshData} shows the relevant mesh and grading parameters for all tested conditions.

A point source was placed on the y-axis close to the surface of the sphere $\bx^*=[0, 101, 0]\,\si{\milli\meter}$. The EQA grid was used as loudspeaker grid.

\begin{table}[!h]
    \centering
    \begin{tabularx}{\linewidth}{Y Y Y Y Y Y Y}
        \toprule
        & \multicolumn{2}{c}{$\hat{\ell}$ (target)} & \multicolumn{3}{c}{$\ell$ (actual)} & \\
        \cmidrule(lr){2-3} \cmidrule(lr){4-6}
        & $\hat{\ell}_{min}$ & $\hat{\ell}_{max}$ & $\ell_{min}$ & $\ell_{max}$ & $\ell_{avg}$ & $\#\mathcal{F}$ \\
        \midrule
        \multirow{6}{*}{UNI} &  \multicolumn{2}{c}{2}  & 1.4 & 2.1 & 1.9 & 81920\\
          &  \multicolumn{2}{c}{3}  & 1.7 & 4.4 & 2.9 & 35490\\
          &  \multicolumn{2}{c}{4}  & 2.5 & 5.3 & 3.8 & 20138\\
          &  \multicolumn{2}{c}{5}  & 3.0 & 6.7 & 4.9 & 12358\\
          &  \multicolumn{2}{c}{10}  & 5.5 & 13.3 & 9.4 & 3278\\
          &  \multicolumn{2}{c}{20}  & 12.1 & 26.7 & 18.7 & 828\\
        \cmidrule(lr){1-7}
        \multirow{5}{*}{POW1} & 2 & 5.1 & 1.4 & 5.2 & 2.8 & 35798\\
          & 2 & 7.5 & 1.4 & 7.7 & 3.6 & 20002\\
          & 2 & 10.5 & 1.5 & 10.5 & 4.6 & 12174\\
          & 2 & 16 & 1.6 & 15.0 & 6.1 & 6432\\
          & 2 & 50 & 1.8 & 42.8 & 14.0 & 994\\
                 \cmidrule(lr){1-7}
        \multirow{5}{*}{POW2} & 2 & 7.2 & 1.2 & 7.3 & 2.7 & 35498\\
          & 2 & 12.5 & 1.4 & 11.9 & 3.3 & 20270\\
          & 2 & 20 & 1.4 & 19.2 & 3.9 & 12550\\
          & 2 & 40 & 1.5 & 34.8 & 4.8 & 6286\\
          & 8 & 60 & 3.9 & 54.1 & 13.9 & 1068\\
        \cmidrule(lr){1-7}
        \multirow{3}{*}{POW4} & 2 & 20 & 1.1 & 18.6 & 2.6 & 29934\\
          & 2 & 40 & 1.1 & 31.1 & 2.8 & 20722\\
          & 4 & 40 & 1.8 & 35.5 & 4.9 & 8160\\
        \cmidrule(lr){1-7}
        \multirow{3}{*}{COS2} & 2 & 5 & 1.1 & 5.2 & 2.8 & 35024\\
          & 2 & 11 & 1.0 & 11.4 & 4.1 & 13474\\
          & 2 & 20 & 1.5 & 20.1 & 8.8 & 2714\\
        \cmidrule(lr){1-7}
        \multirow{4}{*}{COS4} & 2 & 8 & 1.2 & 8.2 & 4.2 & 14420\\
          & 2 & 15 & 1.2 & 15.3 & 6.0 & 5924\\
          & 2 & 20 & 1.0 & 19.7 & 6.8 & 4088\\
        \bottomrule
    \end{tabularx}
    \caption{The SPH objects: Mesh statistics for uniform (UNI) and graded (POW$\alpha$ and COS$\alpha$) meshes. $l_{min}$, $l_{max}$, and $l_{avg}$ are the minimum, maximum, and average edge length in millimeter, respectively. $\#\mathcal{
F}$ denotes the number of elements in a mesh.}
    \label{tab:sphereMeshData}
\end{table}

The relative numerical error (see Eq.~\ref{equ:NumericalError}) was calculated. The reference was an analytically derived HRTF set of the SPH object. To this end, the sound field of a point source scattered by a rigid sphere \citep{beranek_acoustics:_2012} was calculated by:
\begin{equation}
    p_{ref}(\bx)=p_{inc}(\bx)+p_{scat}(\bx),
    \label{equ:pTotal}
\end{equation}
\noindent where $p_{inc}(\bx)=p_0\frac{e^{-\rmi k|\bx-\bx^*|}}{4\pi|\bx-\bx^*|}$ is the incoming sound field,
$p_0$ is the source strength, $k$ is the wavenumber, and $|\bx-\bx^*|$ is the distance between the loudspeaker position and the source. The scattered sound field $p_{scat}$ caused by the sphere is then:
\begin{equation}
\begin{split}
    p_{scat}(\bx)=\frac{\rmi k p_0}{4\pi}\sum_{n=0}^{\infty}
    &(2n+1)h_n^{(2)}(k | \bx^* |)\\
    &\frac{{j_n^\prime}(kR)}{{h_n^\prime}^{(2)}(kR)}h_n^{(2)}(k | \bx |)P_n ( \cos \beta ),
    \label{equ:pScattered}
\end{split}
\end{equation}
\noindent where $j_n$ and $h_n^{(2)}$ are the spherical Bessel and (second order) Hankel functions, respectively, and $j_n^\prime $ and ${h_n^\prime}^{(2)}$ are their derivatives.
$P_n$ are the Legendre polynomials of order $n$. $\cos(\beta)=(\bx^T \bx^*)/(|\bx| |\bx^*|)$ represents the angle between the source and the loudspeaker position. Eq.~\ref{equ:pScattered} was evaluated up to the order of $n=125$.

The relative numerical errors were calculated for the $L^{\mathrm{2}}$-norm (Eq. \ref{equ:NormSphere2}) and the $L^{\mathrm{\infty}}$-norm (Eq.~\ref{equ:NormSphereInf}), both for all, ipsilateral, and contralateral loudspeaker positions. Further, the computation time averaged across calculated frequencies was recorded.

\begin{figure}[!h]
\centering
    \large\rotatebox{90}{\hspace{4.0cm}$e_{L^2}^{\Omega_e}$ in \%}\normalsize
    \includegraphics[width=0.55\textwidth]{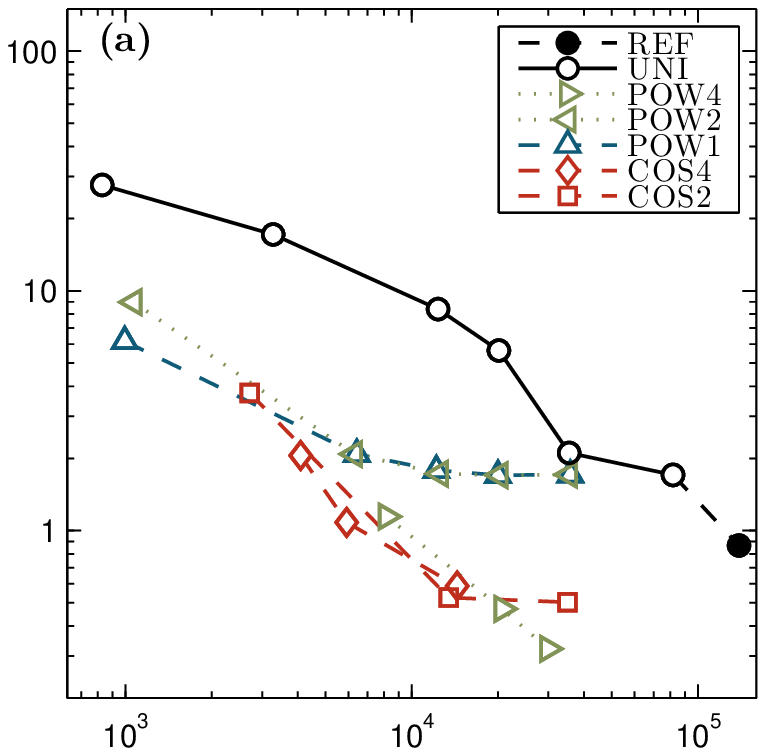}\\
    \vspace{0.2cm}
    \large\rotatebox{90}{\hspace{2.0cm}$e_{L^2}^{\Omega_e}$ in \%}\normalsize
    \hspace{0.2cm}
    \includegraphics[height=0.3\textwidth]{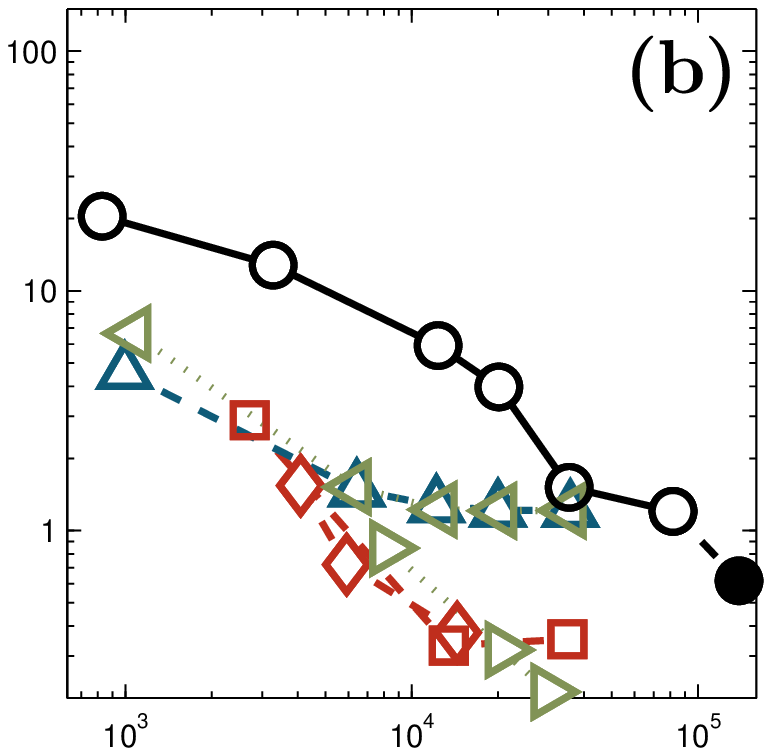}
    \includegraphics[height=0.3\textwidth]{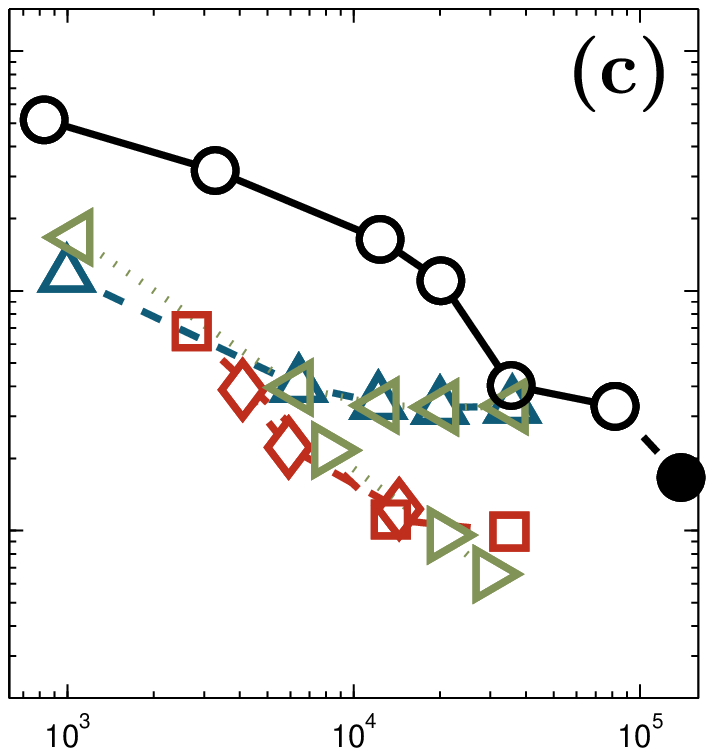}

    \vspace{0.2cm}
    \large \hspace{1.2cm}Number of elements $\#\mathcal{F}$ \normalsize
    \caption{The SPH object and analytical reference: Relative numerical errors for \textbf{(a)} all nodes, \textbf{(b)} ipsilateral nodes, and \textbf{(c)} contralateral nodes of the EQA loudspeaker grid. The reference for the error was the \textit{analytically} derived HRTF set following Eq.~\ref{equ:pTotal}.}
\label{fig:sphereError}
\end{figure}

Fig.~\ref{fig:sphereError} shows the relative numerical errors $e^{\Omega_e}_{L^2}$ as functions of $\#\mathcal{F}$ calculated for various mesh conditions. $e^{\Omega_e}_{L^2}$ considered all (Fig.~\ref{fig:sphereError}a), ipsilateral (b), and contralateral (c) loudspeaker positions. Note that the filled circle highlights the error obtained for the high-resolution mesh, a condition usually thought to provide very accurate results and evaluated as the reference in the following section.

For the uniform meshes (the UNI condition), the error increased monotonically with decreasing $\#\mathcal{F}$. Generally, for the non-uniform meshes (conditions POW$\alpha$ and COS$\alpha$), the error also increased with decreasing $\#\mathcal{F}$. However, in conditions POW1 and POW2, the errors did not decreased much for $\#\mathcal{F}$ beyond \num{8000}. They seem to converge at the level of the error obtained for the UNI condition with an AEL of \SI{2}{\milli\meter}. In conditions COS2, COS4, and POW4, the errors decreased further even for $\#\mathcal{F}$ beyond 8000, showing a similar decay rate as the errors for the UNI condition. These grading functions showed the most promising effect of mesh grading, for which as it seems, $\#\mathcal{F}$ can be reduced by factor of ten without raising the errors significantly.

In conditions COS2, COS4, and POW4, for $\#\mathcal{F}$s beyond \num{10000}, the errors were even smaller than those for both the \SI{2}{\milli\meter} uniform and the high-resolution meshes. This result might appear intriguing because the average resolution of the graded meshes was much smaller than the resolution of those two uniform meshes. It seems like having more elements is not always of advantage, potentially resulting from a larger numerical error when dealing with more elements in the numerical procedures of the BEM (e.g., summing over many elements of an array). In the context of our evaluation, this finding provides evidence that HRTFs based on graded meshes are able to approximate the exact analytical solution at least as good as HRTFs based on uniform meshes.

The errors for the ipsilateral loudspeakers only (Fig.~\ref{fig:sphereError}b) were smaller than those for contralateral ones (Fig.~\ref{fig:sphereError}c). This was expected for the graded meshes, which poorer accuracy is an intrinsic property for the contralateral directions. Interestingly, even in the conditions with uniform meshes, the errors for the ipsilateral loudspeakers were smaller than those for contralateral ones, which might appear surprising given similar element sizes at the two lateral sides in uniform meshes. This finding, however, shows that the elements on the ipsilateral side affect HRTFs much more than those on the contralateral side, supporting our assumption on loosening requirements for the geometric accuracy at the contralateral side.

The pattern of the relative errors based on the $L^{\mathrm{\infty}}$-norm was similar to that based on the $L^{\mathrm{2}}$-norm (thus figures are not shown). Averaged across the tested conditions, $e_{L^{\infty}}^{\Omega_e}$s were only 1.51 times larger than $e_{L^{2}}^{\Omega_e}$s. The similarity between these two errors indicates that 1) none of the loudspeaker positions yielded an error larger than 1.51 of the average error, showing no evidence for an outlier, and 2) the same conclusions can be drawn from both norms.

\subsection{The SPH object and numerical reference}
\label{sub:Experiment12}

\begin{figure}[!h]
\centering
\large\rotatebox{90}{\hspace{4.0cm}$e_{L^2}^{\Omega_e}$ in \%}\normalsize
\includegraphics[height=0.35\textheight]{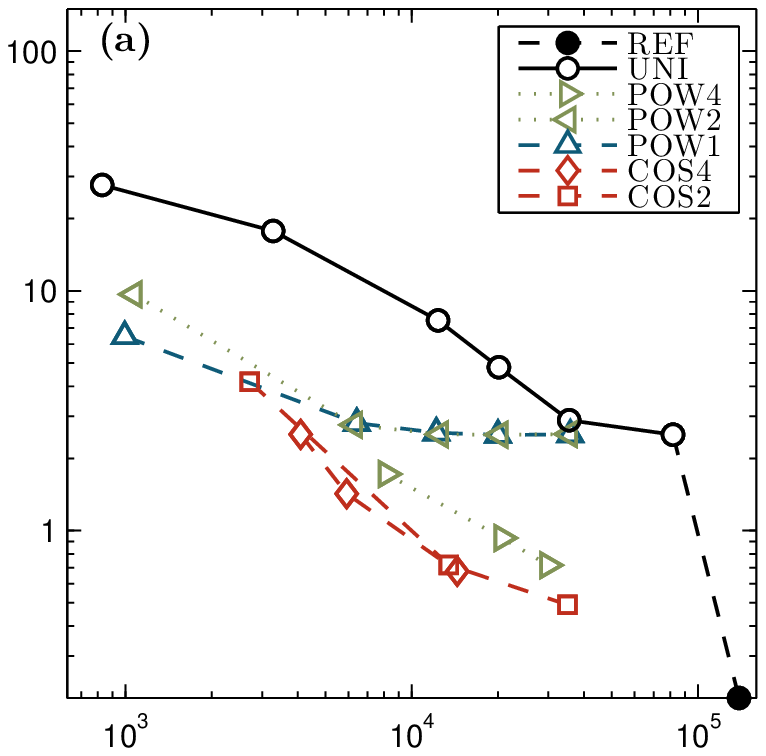}\\
\vspace{0.2cm}
\large\rotatebox{90}{\hspace{2.0cm}$e_{L^2}^{\Omega_e}$ in \%}\normalsize
\hspace{0.23cm}
\includegraphics[height=0.2\textheight]{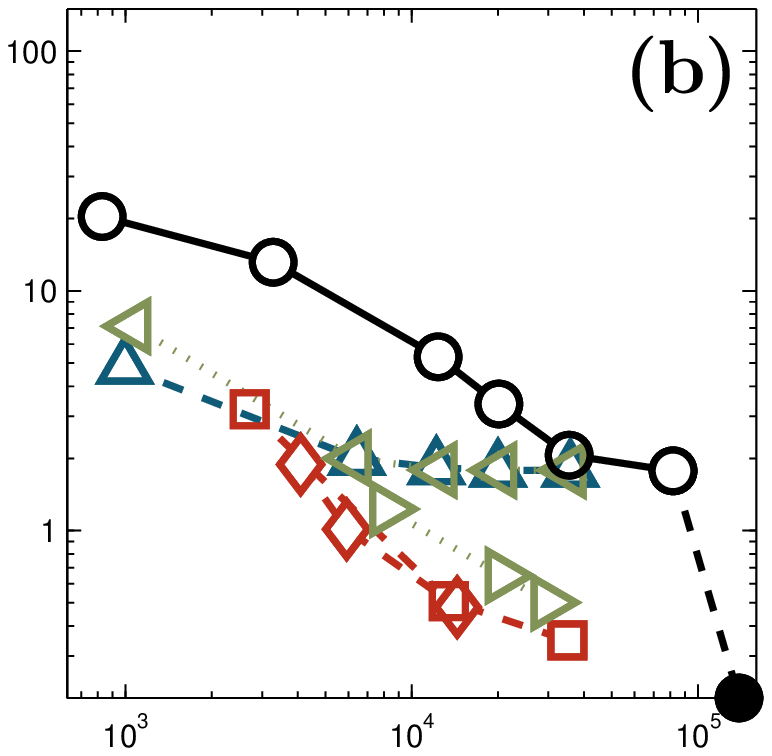}
\includegraphics[height=0.2\textheight]{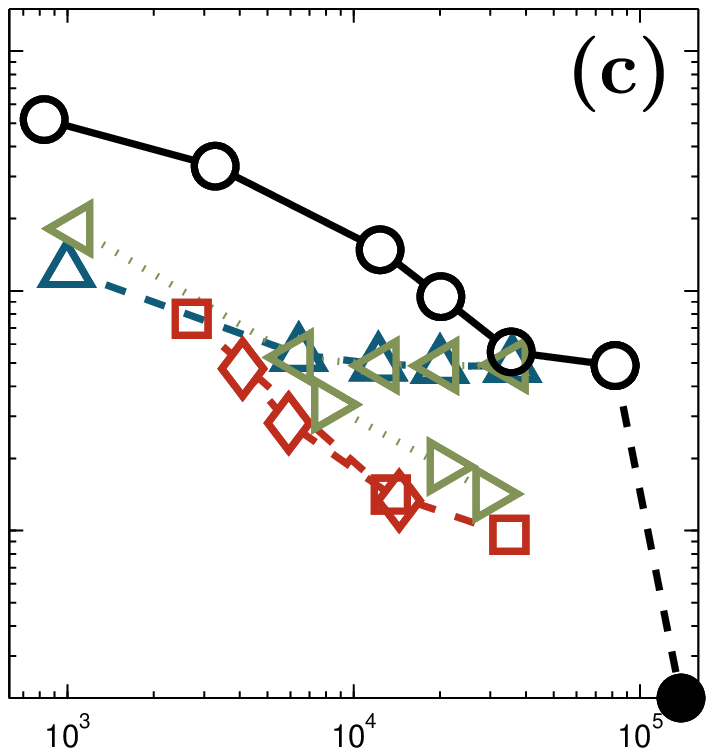}

\vspace{0.2cm}
\large \hspace{1.2cm}Number of elements $\#\mathcal{F}$\normalsize

\caption{The SPH object and numerical reference: Relative numerical errors for \textbf{(a)} all nodes, \textbf{(b)} ipsilateral nodes, and \textbf{(c)} contralateral nodes of the EQA loudspeaker grid. The reference for the errors was the HRTF set \textit{numerically} calculated for the high-resolution mesh.}
\label{fig:sphereErrorRef}
\end{figure}

For the SAP and HUM objects, an analytical derivation of the HRTFs is not feasible and an HRTF set numerically calculated from a high-resolution mesh must be used as a reference. Thus, we investigated the appropriateness of the high-resolution mesh as a basis for reference HRTFs. To this end, the relative errors from the previous section were recalculated with an other reference, namely, the HRTF set calculated for a high-resolution mesh of the SPH object. Further, the gain in computation time from reducing $\#\mathcal{F}$ in the calculations was evaluated.

\begin{figure}[!h]
\centering
\begin{minipage}{8.4cm}
    \flushright
    \large\rotatebox{90}{\hspace{1.8cm}Rel. computation time in $\%$}\normalsize
    \includegraphics[width=7.9cm]{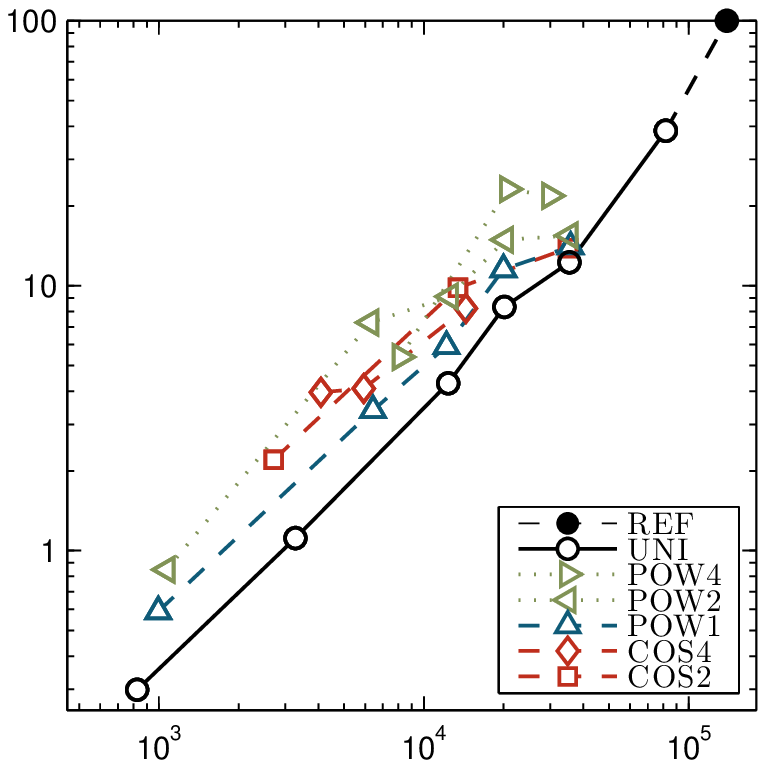}
\end{minipage}

    \vspace{0.2cm}
    \large \hspace{1.5cm} Number of elements $\#\mathcal{F}$ \normalsize

    \caption{Relative computation time to calculate $\phi(\bx)$ in percent. The reference (100\%) was the recorded computation time required for the high-resolution mesh.}
\label{fig:sphereComputationalEffort}
\end{figure}

Fig.~\ref{fig:sphereErrorRef} shows the relative numerical error $e^{\Omega_e}_{L^2}$. Note that the filled circle still highlights the error obtained for the high-resolution mesh, which is zero percent now, and is thus shown on the abscissa. The errors seem to follow similar patterns as those with analytically derived HRTFs as reference (compare Fig.~\ref{fig:sphereError}). Pearson’s correlation coefficient between $e^{\Omega_e}_{L^2}$ calculated with the analytically derived HRTFs and $e^{\Omega_e}_{L^2}$ calculated with the high-resolution mesh HRTFs was \num{0,99}. Both the patterns and the correlation coefficient indicate that the high-resolution mesh adequately represents the reference for further investigations with more complex geometric objects for which an analytical solution is not feasible.

As for $e^{\Omega_e}_{L^{\infty}}$, the errors followed similar patterns as for $e^{\Omega_e}_{L^2}$. When averaged across the tested conditions, $e_{L^{\infty}}^{\Omega_e}$ was only 1.43 times larger than $e_{L^{2}}^{\Omega_e}$, showing no evidence for an outlier across frequencies and loudspeaker positions.

Fig.~\ref{fig:sphereComputationalEffort} shows the relative computation time required to calculate $\phi(\bx)$ as a function of $\#\mathcal{F}$. The reference (i.e., \SI{100}{\percent}) was the computation time for the high-resolution mesh. As expected, the computation time increased with $\#\mathcal{F}$. On average across the tested $\#\mathcal{F}$s, the computation time was slightly larger for graded meshes than for uniform meshes. This might be a side-effect of our FMM implementation, which was optimized for uniform cluster sizes and thus uniform meshes. Compared to the UNI conditions, the POW1 condition showed the smallest increase in computation time, while the conditions POW2 and POW4 showed the largest increase (up to a factor of two) in the required computation time.

An interesting aspect is also the memory required for the calculations. The RAM consumption (for calculating $\phi(\bx)$ at \SI{20}{\kilo\hertz}) was \SI{7,99}{\giga\byte} for the high-resolution mesh. For graded meshes with approximately \num{13000} elements, the RAM consumption was \SI{1,41}{\giga\byte} for the UNI mesh, \SI{1,47}{\giga\byte} for the POW1 mesh, \SI{1,52}{\giga\byte} for the POW2 mesh, and \SI{1,74}{\giga\byte} for the COS2 mesh. Thus, mesh grading not only reduces the computation time, but also loosens the requirements on memory.

In summary, the three grading functions COS2, COS4, and POW4 showed the most potential in terms of smallest numerical errors in HRTF calculations. The POW1 grading showed the most potential in terms of smallest time required to calculate HRTFs.

\subsection{The SAP object}
\label{sub:Experiment2}

\begin{table}[!h]
    \centering
    \begin{tabularx}{\linewidth}{Y Y Y Y Y Y Y}
        \toprule
        & \multicolumn{2}{c}{$\hat{\ell}$ (target)} & \multicolumn{3}{c}{$\ell$ (actual)} & \\
        \cmidrule(lr){2-3} \cmidrule(lr){4-6}
        & $\hat{\ell}_{min}$ & $\hat{\ell}_{max}$ & $\ell_{min}$ & $\ell_{max}$ & $\ell_{avg}$ & $\#\mathcal{F}$ \\
        \midrule
        \multirow{5}{*}{UNI} &  \multicolumn{2}{c}{2}  & 1.1 & 2.8 & 2.0 & 77984\\
          &  \multicolumn{2}{c}{3}  & 1.8 & 4.4 & 2.9 & 35034\\
          &  \multicolumn{2}{c}{5}  & 3.0 & 6.6 & 4.8 & 12962\\
          &  \multicolumn{2}{c}{7}  & 3.9 & 9.9 & 7.0 & 5992\\
          &  \multicolumn{2}{c}{10}  & 6.2 & 13.4 & 10.0 & 2928\\
        \cmidrule(lr){1-7}
        \multirow{3}{*}{POW1} & 2 & 5 & 1.0 & 5.2 & 2.8 & 35162\\
          & 2 & 10 & 1.2 & 10.1 & 4.4 & 13286\\
          & 2 & 20 & 1.3 & 18.9 & 6.8 & 4820\\
        \cmidrule(lr){1-7}
        \multirow{3}{*}{POW4} & 2 & 20 & 1.0 & 18.7 & 2.5 & 32738\\
          & 2 & 40 & 1.0 & 35.9 & 2.6 & 23330\\
          & 3 & 40 & 1.5 & 31.9 & 3.8 & 12968\\
        \cmidrule(lr){1-7}
        \multirow{3}{*}{COS2} & 1 & 12 & 0.5 & 12.3 & 2.4 & 26416\\
          & 1 & 15 & 0.5 & 15.6 & 4.6 & 7494\\
          & 2 & 15 & 1.0 & 15.4 & 6.7 & 4730\\
        \bottomrule
    \end{tabularx}
    \caption{The SAP object: Mesh statistics for uniform (UNI) and graded (POW1, POW4, and COS2) meshes. Other details as in Tab.~\ref{tab:sphereMeshData}.}
    \label{tab:spherePinnaMeshData}
\end{table}

The SAP object (SPH object with a generic pinna) was used to investigate the effect of an additional pinna on the mesh grading with a good relation to the results obtained for the SPH object. In order to create the SAP object, the default pinna from MakeHuman\footnote{MakeHuman: version 1.0.0,  available from \url{http://www.makehuman.org} (date last viewed: January 31, 2016)} was stitched onto the left side, i.e., the positive y dimension, of the SPH object with the entrance of the ear canal at position $\bx=[0, 100, 0]\,\si{\milli\meter}$. The high-resolution mesh consisted of \num{140847}~elements with an AEL of approximately \SI{1,4}{\milli\meter}. The high-resolution mesh was further re-meshed in order to obtain uniform meshes with AELs ranging from \num{2} to \SI{10}{\milli\meter}. Also, the a-priori mesh grading was applied to the high-resolution mesh. The grading functions POW1, POS4, and COS2 were used. The COS4 and POW2 functions were not used because their previous results were very similar to those for COS2 and POW1, respectively. Table~\ref{tab:sphereMeshData} shows the relevant mesh parameters for the tested grading functions and conditions.

For each mesh, a point source was placed at $x^*=[0, 101, 0]\,\si{\milli\meter}$ and HRTFs were calculated using the EQA loudspeaker grid. The relative numerical errors with the high-resolution mesh HRTFs as reference were calculated.

\begin{figure}[!h]
    \centering
    \includegraphics[width=0.7\linewidth]{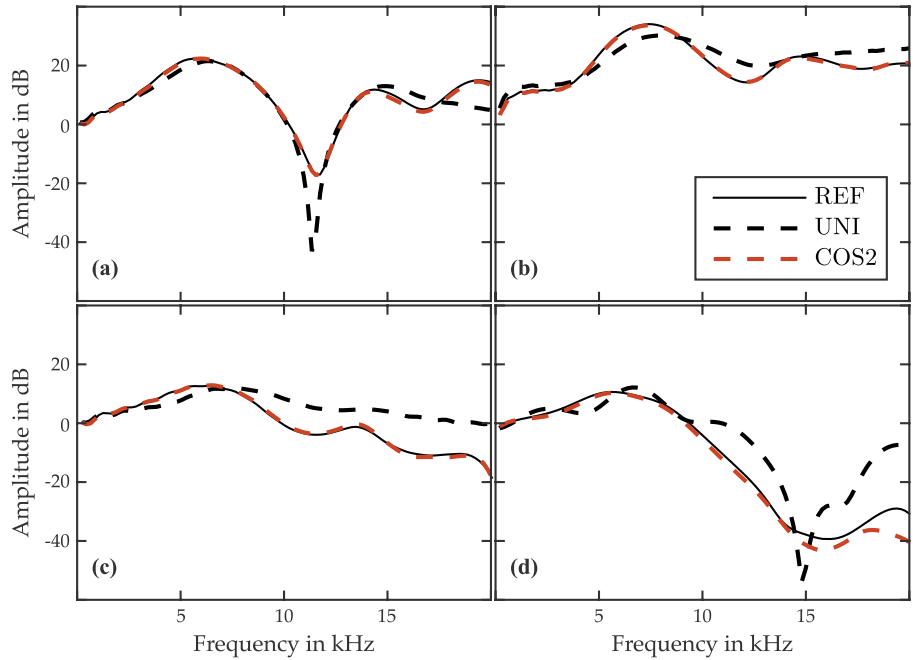}
    \caption{HRTF amplitude spectra calculated for the SPH object represented by the high-resolution mesh (REF, \num{108940}~elements), uniform mesh (UNI, \num{12962}~elements), and graded mesh (COS2, \num{7494}~elements). The loudspeaker positions were \textbf{(a)} \ang{0}, \textbf{(b)} \ang{90}, \textbf{(c)} \ang{180}, and \textbf{(d)} \ang{270} (contralateral) in the horizontal plane.}
 \label{fig:spherePinnaSolution}
\end{figure}

Fig.~\ref{fig:spherePinnaSolution} shows exemplary HRTFs of a uniform mesh (UNI, AEL of \SI{5}{\milli\meter}) and a non-uniform mesh (COS2, minimum and maximum target edge length of \num{1} and \SI{15}{\milli\meter}, respectively) calculated for a low $\#\mathcal{F}$ at four loudspeaker positions in the horizontal plane (front, left, rear, and right). The HRTFs calculated for the high-resolution mesh (REF) served as a reference in that comparison. While HRTFs of the uniform mesh differ in general at frequencies above \SI{5}{\kilo\hertz}, HRTFs of the graded mesh seem to correspond very well with the reference HRTFs. Only for the contralateral position, differences arise above \SI{15}{\kilo\hertz} but still these differences are in the range of \SI{2}{\dB} on average. For that position, HRTFs of the uniform mesh, show obvious differences, with an additional notch at \SI{15}{\kilo\hertz}. Note that in this comparison, a rather low $\#\mathcal{F}$ (about ten thousand elements) was used in order to 1) underline the limits of uniform meshes for the numerical HRTF calculation, and 2) evaluate the spectral details in HRTFs based on graded meshes.

\begin{figure}[!h]
\centering
\begin{minipage}{8.4cm}
    \flushright
    \large\rotatebox{90}{\hspace{3.5cm}$e_{L^2}^{\Omega_e}$ in \%}\normalsize
    \includegraphics[height=7.0cm]{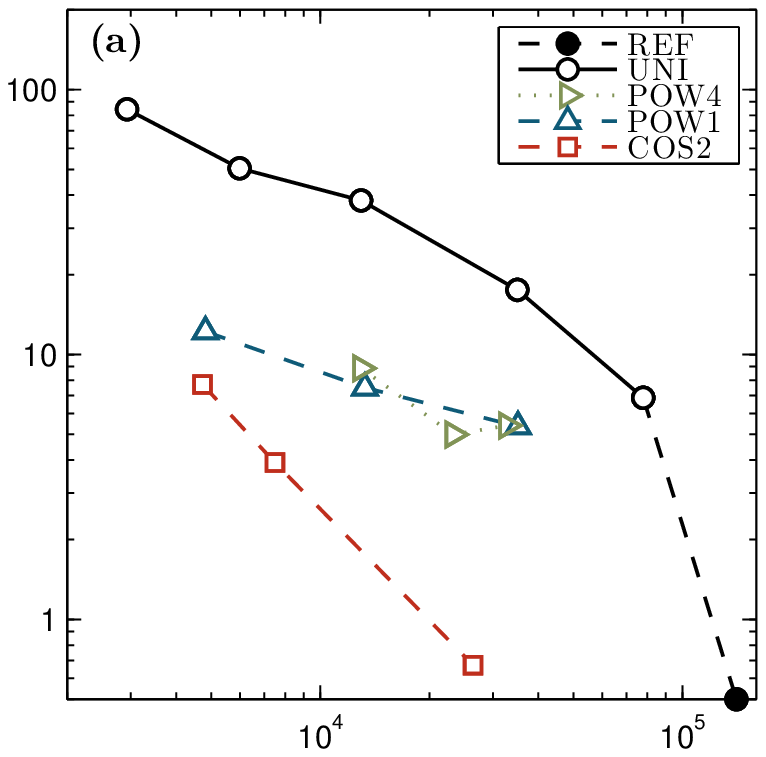}\\
    \vspace{0.2cm}
    \large\rotatebox{90}{\hspace{1.4cm}$e_{L^2}^{\Omega_e}$ in \%}\normalsize
    \hspace{0.227cm}
    \includegraphics[height=3.5cm]{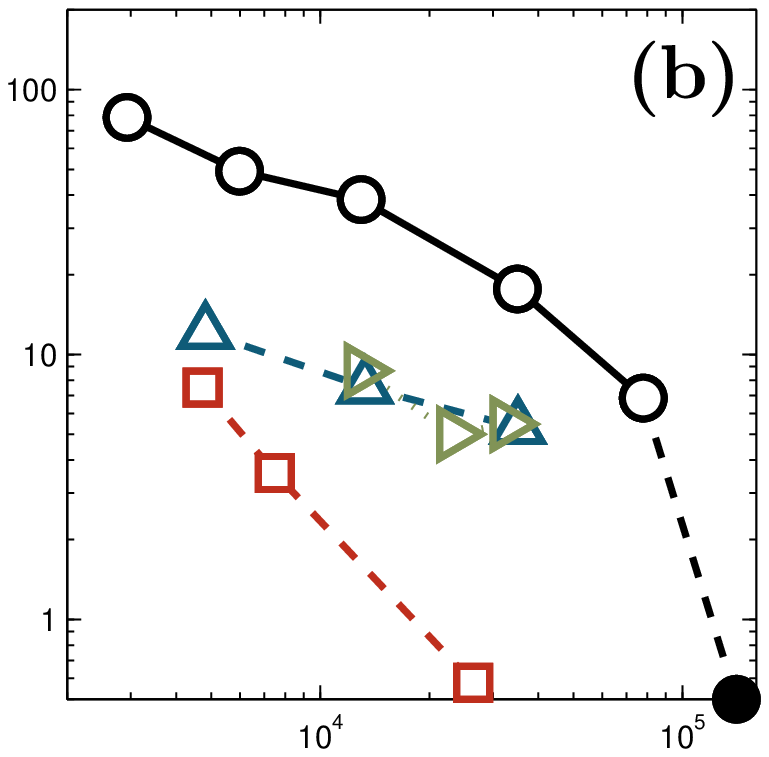}
    \includegraphics[height=3.5cm]{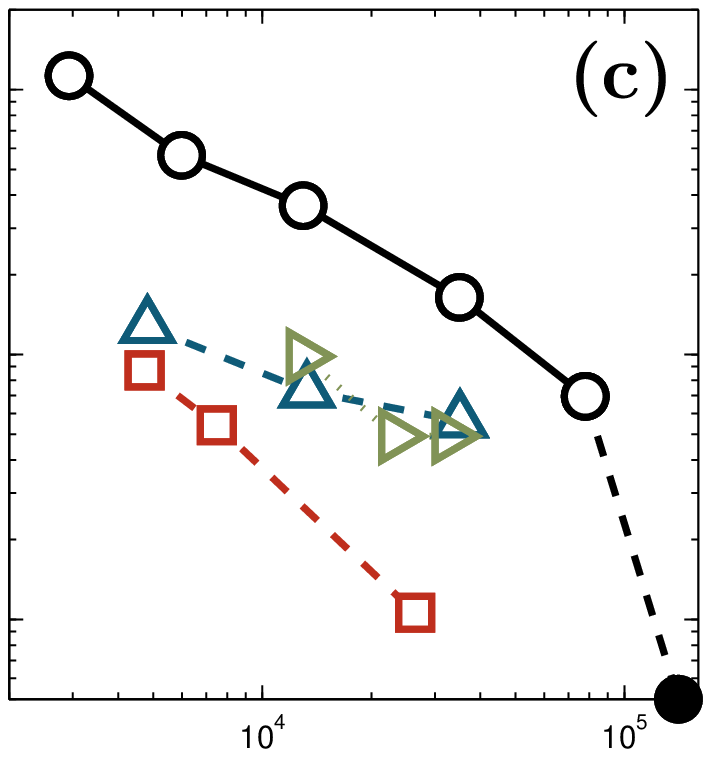}
\end{minipage}

    \vspace{0.2cm}
    \large \hspace{1.2cm}Number of elements $\#\mathcal{F}$\normalsize

    \caption{The SAP object: Relative numerical errors for \textbf{(a)} all nodes, \textbf{(b)} ipsilateral nodes, and \textbf{(c)} contralateral nodes of the EQA grid. The reference for the errors was the HRTF set calculated for the high-resolution mesh.}
\label{fig:spherePinnaError}
\end{figure}

Fig.~\ref{fig:spherePinnaError} shows $e_{L^2}^{\Omega_e}$ calculated for all (a), ipsilateral (b), and contralateral (c) loudspeakers. Similar to results for the SPH object, the errors increased with decreasing $\#\mathcal{F}$. Also, for similar $\#\mathcal{F}$, POW1 and POW4 conditions yielded similarly smaller errors than the UNI condition. The COS2 condition yielded even smaller errors. The POW$\alpha$ conditions seem to require three times more and the UNI condition seems to require ten times more elements than the COS2 condition to achieve a similar relative error in the calculated HRTFs.

Averaged across all conditions, the $e_{L^{\infty}}^{\Omega_e}$ was \num{1,09} times larger than $e_{L^2}^{\Omega_e}$ showing a very small variance and no evidence for an outlier across the loudspeaker positions.

In summary, the meshes graded with the POW4 function yielded errors similar to those of the meshes graded with the POW1 function. The COS2 grading function showed the most promising results. For this function, the HRTF spectra did not show any problematic issues. And even for the coarsest COS2 graded mesh (\num{4730}~elements), the error was smaller than that or HRTFs calculated for the finest uniform mesh (\num{77984}~elements, AEL of \SI{2}{\milli\meter}). This is interesting because the minimum edge length of that COS2 graded mesh was similar to the edge length of that uniform mesh (see $\ell_{min}$ in Table~\ref{tab:spherePinnaMeshData}), implying that while both meshes had a similar resolution at the ipsilateral side, the graded mesh had a much coarser resolution at the contralateral side. This supports again our assumption that a fine mesh resolution at the ipsilateral pinna is important, whereas for the rest of the mesh, a coarser resolution can still yield acceptable accuracy.

\subsection{The HUM objects}
\label{sub:Experiment3}
The HUM objects were geometric models of the head and pinnae of three actual listeners (NH5, NH130, and NH131). These high-resolution meshes originate from scanning the head with a laser scanner and scanning silicone impressions of listener's pinnae in a high-energy industrial computer-tomography scanner \citep{reichinger_evaluation_2013}. The meshes consisted of \numlist{111362;111422;107692} elements for NH5, NH130, and NH131, respectively, and are described in more details in \citet{ziegelwanger_numerical_2015}.

For each listener, the high-resolution mesh was re-meshed in order to obtain uniform meshes with AELs ranging from \SIrange{2}{10}{\milli\meter}. Graded meshes were created by applying the a-priori mesh grading to the high-resolution meshes\footnote{Note that for each ear separate meshes were created, yielding two meshes per subject and condition.}. Based on the results from Sec.~\ref{sub:Experiment12} and \ref{sub:Experiment2}, the grading functions POW1 and COS2 were considered only. Table~\ref{tab:subjectsMeshData} shows relevant mesh parameters for the tested grading functions and conditions. For each mesh, HRTFs were calculated for both grids (ARI and EQA). Then, the relative numerical errors with the high-resolution mesh HRTFs as reference were calculated.

\begin{table}[!h]
    \centering
    \begin{tabularx}{\linewidth}{Y Y Y Y Y Y Y}
        \toprule
        & \multicolumn{2}{c}{$\hat{\ell}$ (target)} & \multicolumn{3}{c}{$\ell$ (actual)} & \\
        \cmidrule(lr){2-3} \cmidrule(lr){4-6}
        & $\hat{\ell}_{min}$ & $\hat{\ell}_{max}$ & $\ell_{min}$ & $\ell_{max}$ & $\ell_{avg}$ & $\#\mathcal{F}$ \\
        \midrule
        \multirow{5}{*}{UNI} &  \multicolumn{2}{c}{2}  & 0.3 & 3.8 & 1.9 & 98009\\
          &  \multicolumn{2}{c}{3}  & 1.5 & 4.3 & 2.8 & 50115\\
          &  \multicolumn{2}{c}{4}  & 2.2 & 5.6 & 3.7 & 28320\\
          &  \multicolumn{2}{c}{5}  & 2.8 & 19.7 & 5.2 & 16630\\
          &  \multicolumn{2}{c}{10}  & 2.9 & 19.7 & 6.3 & 4689\\
        \cmidrule(lr){1-7}
        \multirow{4}{*}{POW1} & 2 & 5 & 1.0 & 5.5 & 3.0 & 39447\\
          & 2 & 8 & 1.1 & 8.4 & 3.6 & 25774\\
          & 2 & 10 & 1.1 & 10.5 & 5.2 & 12407\\
          & 2 & 12 & 1.1 & 12.4 & 5.7 & 9768\\
        \cmidrule(lr){1-7}
        \multirow{3}{*}{COS2} & 1 & 5 & 0.5 & 5.2 & 2.3 & 57526\\
          & 1 & 8 & 0.5 & 8.4 & 2.9 & 30697\\
          & 1 & 15 & 0.5 & 16.5 & 3.7 & 13722\\
        \bottomrule
    \end{tabularx}
    \caption{The HUM objects: Mesh statistics for uniform (UNI) and graded (POW1 and COS2) meshes averaged across subjects and ears. Other details as in Tab.~\ref{tab:sphereMeshData}.}
    \label{tab:subjectsMeshData}
\end{table}

\begin{figure}[!h]
    \centering
    \includegraphics[width=0.9\linewidth]{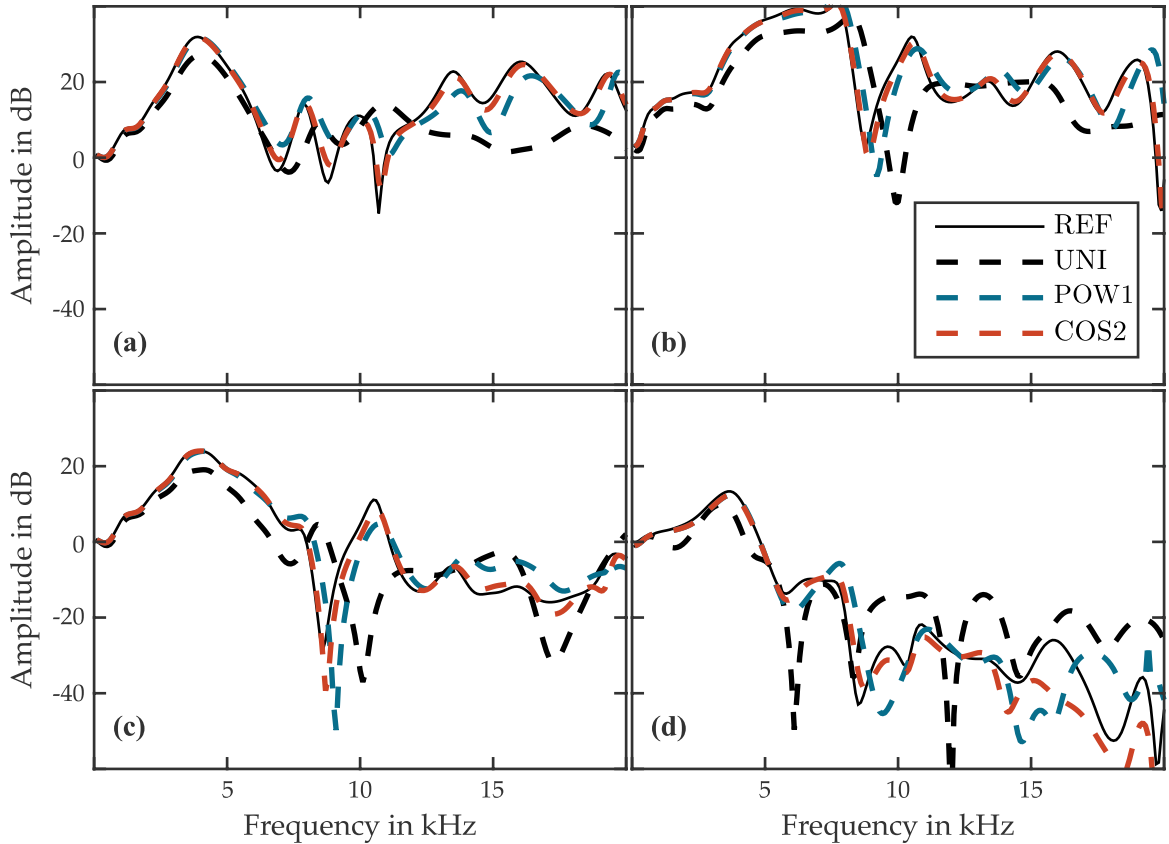}
    \caption{HRTF amplitude spectra calculated for the HUM object of NH5 from the ARI database, represented by the high-resolution mesh (REF, \num{111362}~elements), a uniform mesh (UNI, \num{17023}~elements), a linearly graded mesh (POW1, \num{12041}~elements), and raised cosine graded mesh (COS2, \num{13633}~elements). The loudspeaker positions were \textbf{(a)} \ang{0}, \textbf{(b)} \ang{90}, \textbf{(c)} \ang{180}, and \textbf{(d)} \ang{270} in the horizontal plane.}
\label{fig:subjectsSolution}
\end{figure}

Fig.~\ref{fig:subjectsSolution} shows HRTFs of NH5 calculated for the meshes with between \num{10000} and \num{20000} elements and four loudspeaker positions as in Fig.~\ref{fig:spherePinnaSolution}. The high-resolution mesh HRTFs are shown as a reference. HRTFs calculated for the COS2-graded meshes seem to correspond very well to HRTFs of the high-resolution mesh, except for minor differences at the contralateral loudspeaker position above \SI{15}{\kilo\hertz}. HRTFs calculated for the UNI mesh show much more deviation, e.g., an additional notch at around \SI{10}{\kilo\hertz} arose at the ipsilateral and rear loudspeaker positions.

\begin{figure}[!h]
\centering
\begin{minipage}{8.4cm}
    \flushright
    \large\rotatebox{90}{\hspace{3.5cm}$e_{L^2}^{\Omega_e}$ in \%}\normalsize
    \includegraphics[height=7.5cm]{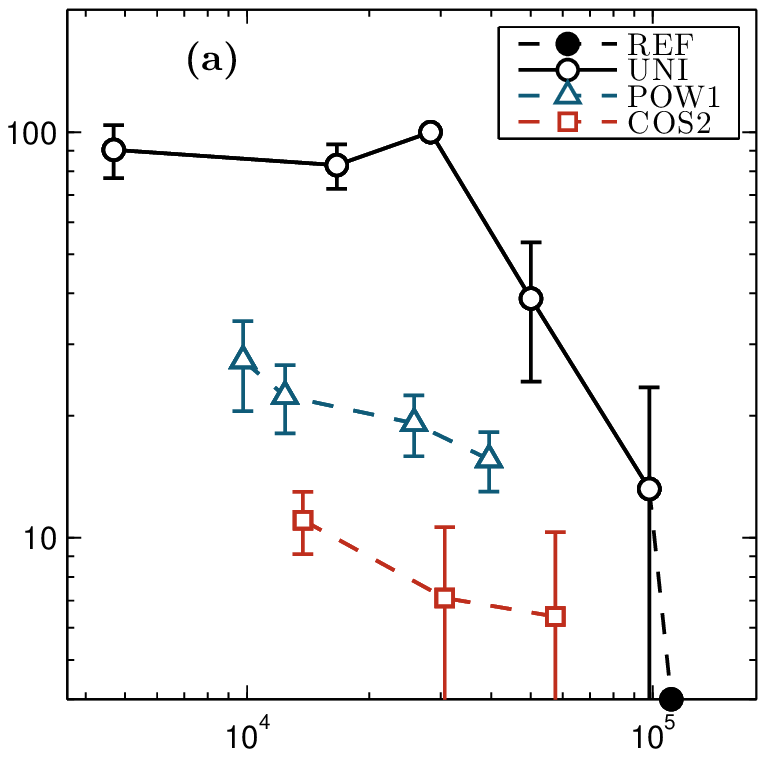}\\
    \vspace{0.2cm}
    \large\rotatebox{90}{\hspace{1.4cm}$e_{L^2}^{\Omega_e}$ in \%}\normalsize
    \hspace{0.227cm}
    \includegraphics[height=3.605cm]{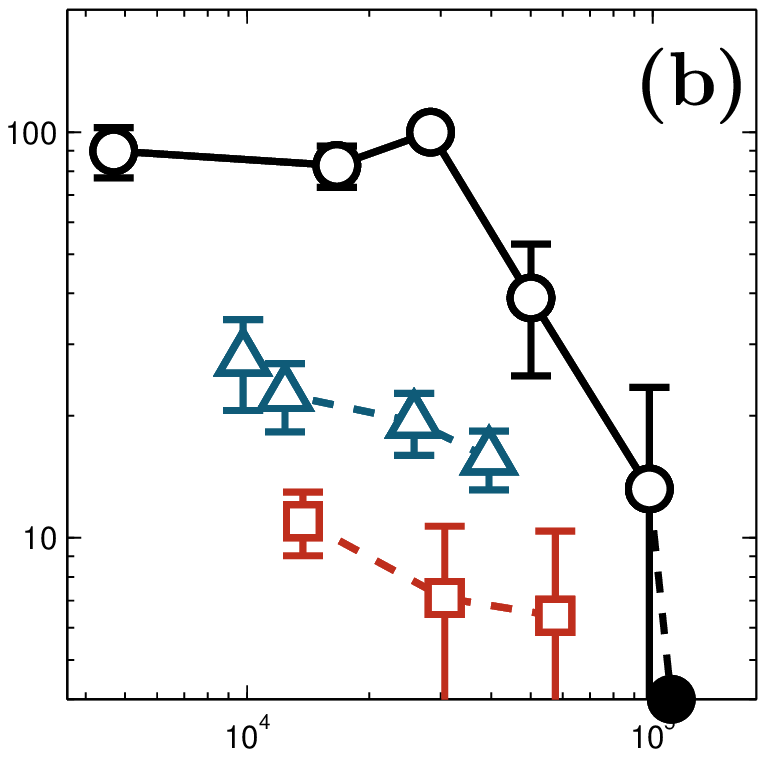}
    \includegraphics[height=3.605cm]{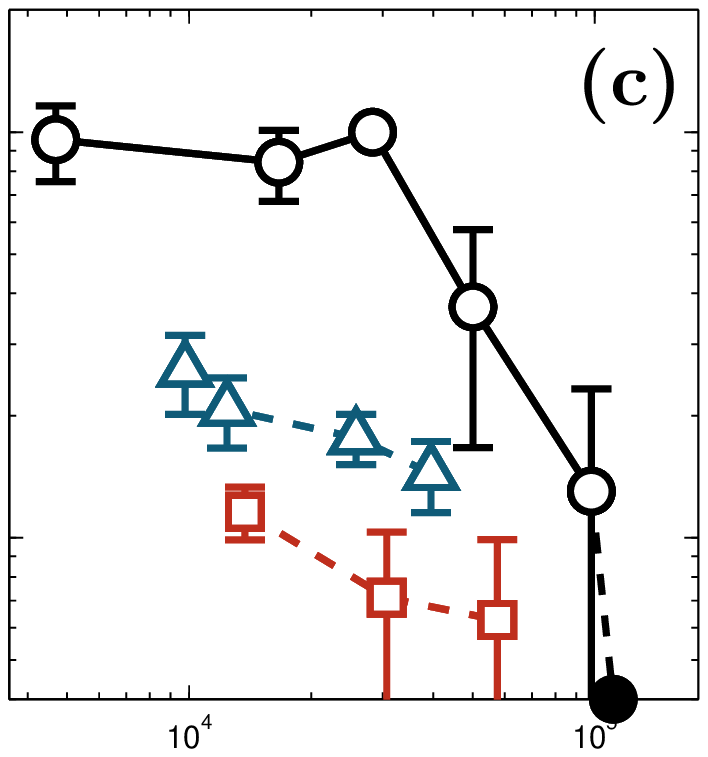}
\end{minipage}

    \vspace{0.2cm}
    \large \hspace{1.2cm}Number of elements $\#\mathcal{F}$\normalsize
    \caption{HUM objects: Relative numerical errors averaged across the listeners and ears for \textbf{(a)} all nodes, \textbf{(b)} ipsilateral, and \textbf{(c)} contralateral positions of the EQA grid. Reference for the error was the HRTF set calculated for the high-resolution mesh of each individual listener and ear.}
\label{fig:subjectsError}
\end{figure}

Fig.~\ref{fig:subjectsError} shows $e_{L^2}^{\Omega_e}$ averaged across subjects and ears and calculated for all (a), ipsilateral (b), and contralateral (c) loudspeakers of the EQA grid. For uniform meshes, the errors increased with decreasing $\#\mathcal{F}$, reaching a maximum at around 100\% for $\#\mathcal{F}$ of approximately \num{28000}, and saturating at this level for less elements. For the POW1-graded meshes, the errors increased with decreasing $\#\mathcal{F}$, but the effect of $\#\mathcal{F}$ was not that large as for uniform meshes and for other objects. For the COS2-graded meshes, the effect of $\#\mathcal{F}$ was more pronounced. Compared to the errors in the UNI condition, the relative error in the COS2 condition was approximately ten times smaller. For the coarsest COS2-graded mesh ($\#\mathcal{F} \approx$ \num{13000}), the error was even slightly smaller than for the finest uniform mesh ($\#\mathcal{F} \approx$ \num{98000}).

For the psychoacoustic-motivated evaluation, acoustically measured HRTFs of the corresponding listeners from the ARI database were used as references\footnote{Available from \url{http://www.kfs.oeaw.ac.at/hrtf} (date last viewed: January 31, 2016)}. These HRTFs were measured in a semi-anechoic chamber by placing in-ear microphones in the listeners' blocked ear canals and applying the multiple exponential sweep method \citep{majdak_multiple_2007} for the system identification of each HRTF direction. These HRTFs were available at positions described by the ARI grid \citep[for more details on the measurement see][]{majdak_3-d_2010,ziegelwanger_numerical_2015}.

First, broadband temporal features were evaluated by means of equivalent head radius derived for an HRTF set (for each ear separately) from the TOA model. Table~\ref{tab:subjectsTOA} shows the equivalent head radii obtained for all conditions, listeners, and ears. For the most of the numerically calculated HRTFs, the equivalent head radii were $\pm$~\SI{2}{\milli\meter} around those of measured HRTFs, showing no evidence for artifacts larger than a few $\mu$\si{\second} in the broadband timing of the calculated HRTFs. The equivalent head radii also did not change much when $\#\mathcal{F}$ was reduced, indicating that the proposed mesh grading did not introduce critical artifacts in the broadband temporal features of calculated HRTFs.

\begin{table}[!h]
    \centering
    \begin{tabularx}{\linewidth}{Y Y Y Y Y Y Y Y Y}
        \toprule
        & & & \multicolumn{6}{c}{Equivalent head radius in mm} \\
        \cmidrule(lr){4-9}
        & \multicolumn{2}{c}{$\hat{\ell}$ (target)} & \multicolumn{2}{c}{NH5} & \multicolumn{2}{c}{NH130} & \multicolumn{2}{c}{NH131} \\
        \cmidrule(lr){2-3} \cmidrule(lr){4-5} \cmidrule(lr){6-7} \cmidrule(lr){8-9}
        & $\hat{\ell}_{min}$ & $\hat{\ell}_{max}$ & l & r & l & r & l & r \\
        \midrule
        \multirow{1}{*}{AC} &  - & -  & 91.8 & 92.4 & 92.4 & 92.0 & 94.5 & 94.8\\
        \cmidrule(lr){1-9}
        \multirow{6}{*}{UNI} &  \multicolumn{2}{c}{1}  & 90.8 & 90.8 & 89.2 & 89.5 & 94.6 & 94.6\\
          &  \multicolumn{2}{c}{2}  & 91.4 & 91.4 & 89.4 & 89.7 & 94.6 & 94.6\\
          &  \multicolumn{2}{c}{3}  & 91.4 & 91.7 & 89.4 & 89.4 & 94.6 & 94.6\\
          &  \multicolumn{2}{c}{4}  & 91.8 & 91.6 & 90.2 & 90.2 & 94.5 & 94.5\\
          &  \multicolumn{2}{c}{5}  & 91.1 & 91.1 & 89.8 & 89.8 & 95.1 & 95.1\\
          &  \multicolumn{2}{c}{10}  & 88.2 & 89.0 & 88.3 & 88.5 & 92.8 & 92.8\\
        \cmidrule(lr){1-9}
        \multirow{4}{*}{POW1} & 2 & 5 & 91.3 & 91.2 & 89.3 & 89.3 & 94.7 & 94.7\\
          & 2 & 8 & 90.9 & 90.9 & 89.5 & 89.5 & 94.7 & 94.7\\
          & 2 & 10 & 91.4 & 91.3 & 89.4 & 89.4 & 94.3 & 94.3\\
          & 2 & 12 & 90.9 & 90.7 & 88.8 & 89.1 & 94.6 & 94.6\\
        \cmidrule(lr){1-9}
        \multirow{3}{*}{COS2} & 1 & 5 & 90.7 & 90.7 & 89.1 & 89.1 & 94.6 & 94.6\\
          & 1 & 8 & 90.7 & 90.6 & 89.0 & 89.2 & 94.3 & 94.3\\
          & 1 & 15 & 90.3 & 90.3 & 88.8 & 88.8 & 94.2 & 94.2\\
        \bottomrule
    \end{tabularx}
    \caption{Equivalent head radii calculated separately for the (\textbf{l})eft and and (\textbf{r})ight ear of NH5, NH130, and NH131 from the ARI database. Acoustically measured HRTFs (AC) and HRTFs numerically calculated for uniform (UNI) and graded meshes (LIN and COS2). Other details as in Tab.~\ref{tab:sphereMeshData}.}
    \label{tab:subjectsTOA}
\end{table}

Second, spectral features were evaluated by means of sound-localization performance predicted for the calculated HRTFs by the sagittal-plane sound-localization model. The performance predictions were calculated for the three listeners. Besides the reference HRTFs (which were the acoustically measured HRTFs), the model requires a parameter called 'uncertainty', representing the ability of a listener to localize sounds \citep{baumgartner_assessment_2013}. We used an uncertainty of \SI{1.9}, corresponding to an average localizer. Four benchmarks were used for the analysis. The first benchmark, 'ACTUAL', was the actual localization performance of \num{14} human listeners \citep{middlebrooks_individual_1999}. Comparison to that benchmark allows to estimate the quality of predictions relative to the actual localization performance. The second benchmark, 'OWN', was the localization performance predicted for \num{177} listeners from three HRTF databases (ARI, LISTEN, CIPIC)\footnote{Available from \url{http://sofaconventions.org} (date last viewed: January 31, 2016)} localizing sound sources with their own ears' HRTFs. Comparison to that benchmark allows to estimate the usual range of predictions when localizing with own ears. The third benchmark, 'OTHER' was the localization performance predicted for our three listeners localizing sound sources with \num{176} others' HRTF sets from the databases. This benchmark represents the result of localizing with non-individual HRTFs. Reaching this error level indicates no need for individual HRTFs in the corresponding condition \citep[for more details see ][]{ziegelwanger_numerical_2015}. The fourth benchmark, 'REF' was the localization performance predicted for the three listeners localizing sources with the HRTFs calculated for the high-resolution mesh. Comparison to this benchmark allows to estimate the ultimate effect of re-meshing on the sound-localization performance.

\begin{figure}[!h]
    \centering
    \includegraphics[width=0.85\linewidth]{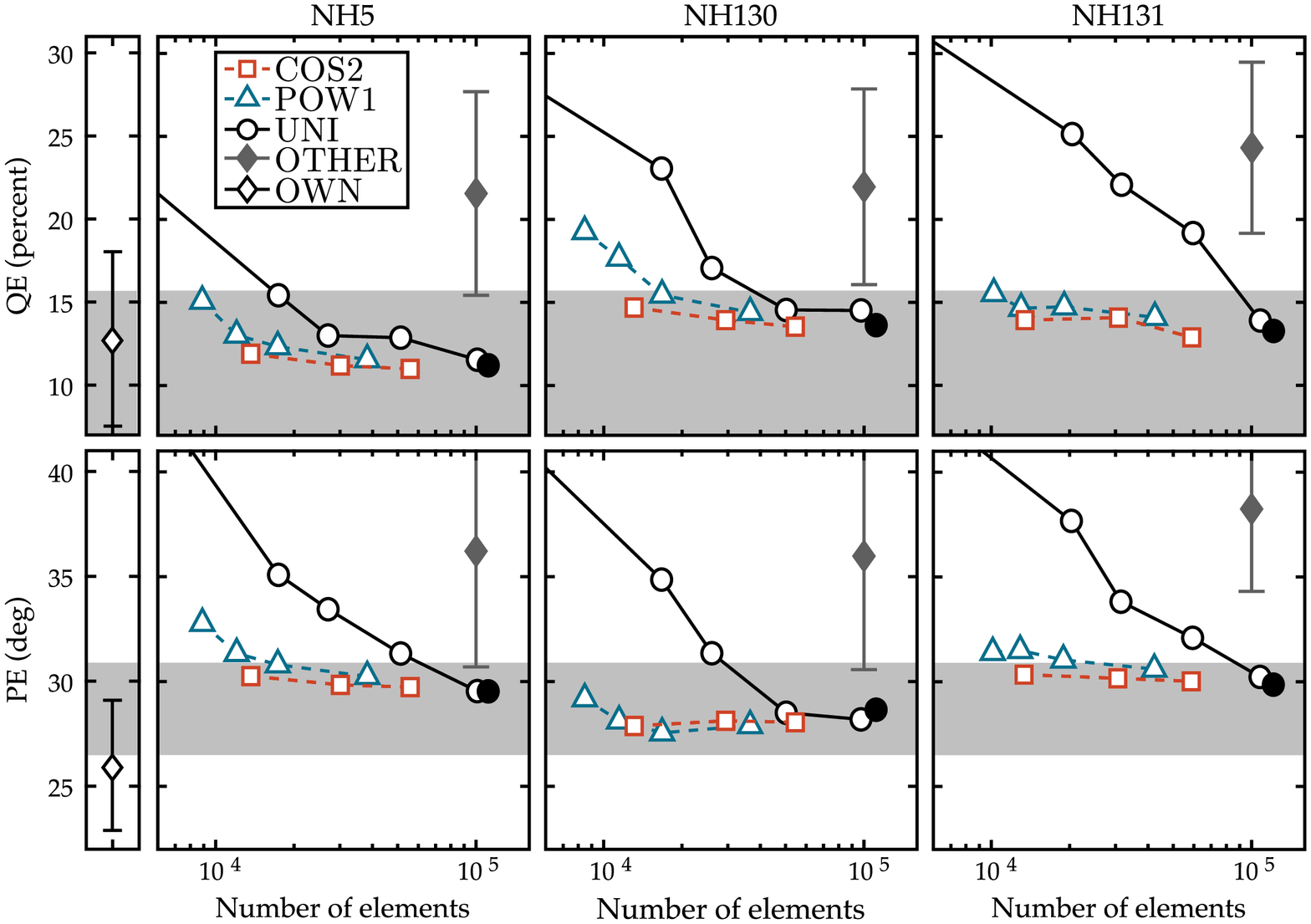}
    \caption{Effect of re-meshing on the sagittal-plane localization performance by means of quadrant error (QE) and polar error (PE). Most-left column: Benchmark 'OWN'. Other columns (NH5, NH130, and NH131): Predictions for the corresponding listener localizing with own-ear HRTFs calculated for the uniform (open circles) and for the non-uniform (triangles) meshes as functions of $\#\mathcal{F}$. Filled circles: Benchmark 'REF'. Filled diamonds: Benchmark 'OTHER'. Grey area: benchmark 'ACTUAL'. See text for details on the benchmarks.}
\label{fig:subjectsPerceptualError}
\end{figure}

Fig.~\ref{fig:subjectsPerceptualError} shows the predicted localization performance as functions of $\#\mathcal{F}$ along with the benchmarks. The benchmarks 'OWN' and 'OTHER' were within and outside the actual listener performance represented by the benchmark 'ACTUAL', respectively, implying a good level of validity for the predictions. For each of the listeners, the benchmark 'REF' was within the range of 'ACTUAL', confirming the general ability of numerically calculated HRTFs to replace acoustically measured HRTFs. For uniform meshes, the performance started at the level of the benchmark 'REF' for the finest mesh and degraded with decreasing $\#\mathcal{F}$. For the coarsest uniform mesh, the predictions were at (or beyond) the level of the benchmark 'OTHER', indicating that these individual HRTFs did not provide any advantages over non-individual HRTFs. For graded meshes, the performance also started at the level of the benchmark 'REF' and degraded with decreasing $\#\mathcal{F}$. In contrast to the uniform meshes, the effect of $\#\mathcal{F}$ was much smaller and the performance was still within the range of benchmark 'OWN' down to meshes with $\#\mathcal{F}$ of approximately \num{13000}.
Compared to uniform meshes, the mesh grading allowed to reduce $\#\mathcal{F}$ by a factor of approximately seven without introducing a degradation of the predicted performance. When compared across the two grading functions, the COS2 grading seems to provide a slightly better performance than the POW1 grading.

In order to roughly quantify the relation between the predicted localization performance and $\#\mathcal{F}$, linear regressions were fit to the predictions.
Similar to the presentation of the abscissa in Fig.~\ref{fig:subjectsPerceptualError}, the fits were performed on $\log_{10} (\#\mathcal{F})$ yielding a PE slope (in $\scriptsize \raisebox{4pt}{$\rm o$} \normalsize  / \scriptsize \raisebox{-0.5pt}{$\log \#\mathcal{F}$} \normalsize$) and a QE slope (in $\scriptsize \raisebox{2pt}{$\%$} \normalsize  / \scriptsize \raisebox{-0.5pt}{$\log \#\mathcal{F}$} \normalsize$) per listener and mesh condition.
The individual slopes and their averages are shown in Table~\ref{tab:subjectsSlopes}. For example, for NH5 and uniform meshes, the PE slope of $-11.2\,\scriptsize \raisebox{4pt}{$\rm o$} \normalsize  / \scriptsize \raisebox{-0.5pt}{$\log \#\mathcal{F}$} \normalsize$ indicates that by increasing $\#\mathcal{F}$ by the factor of ten, the PE decreased by \ang{11,2}.

Our results show that when reducing $\#\mathcal{F}$ in the uniform mesh by factor of ten, the PE and QE errors increased by \ang{10} and \SI{10}{\percent}, respectively. In contrast, for the COS2-graded meshes, the same reduction of elements increased these errors by \ang{0,23} and \SI{1,7}{\percent}, respectively, only.

 \begin{table}[!h]
    \centering
    \setlength{\tabcolsep}{0.5mm}
    \begin{tabularx}{\linewidth}{Y Y Y Y Y Y Y Y Y}
        \toprule
        & \multicolumn{8}{c}{PE ($\scriptsize \raisebox{3pt}{$\rm o$} \normalsize  / \scriptsize \raisebox{-0.5pt}{$\log \#\mathcal{F}$} \normalsize$) and QE ($\scriptsize \raisebox{2pt}{$\%$} \normalsize  / \scriptsize \raisebox{-0.5pt}{$\log \#\mathcal{F}$} \normalsize$) slopes} \\
        \cmidrule(lr){2-9}
        & \multicolumn{2}{c}{NH5} & \multicolumn{2}{c}{NH130} & \multicolumn{2}{c}{NH131} & \multicolumn{2}{c}{Avg.} \\
        \cmidrule(lr){2-3} \cmidrule(lr){4-5} \cmidrule(lr){6-7} \cmidrule(lr){8-9}
        & PE & QE & PE & QE & PE & QE & PE & QE \\
        \midrule
        UNI & -11.2 & -8.2 & -9.8 & -11.0 & -10.0 & -13.2 & -10.2 & -10.5\\
        \cmidrule(lr){1-9}
        POW1 & -3.6 & -4.9 & -1.7 & -7.5 & -1.4 & -1.9 & -1.8 & -5.0\\
        \cmidrule(lr){1-9}
        COS2 & -0.9 & -1.5 & 0.3 & -1.9 & -0.5 & -1.5 & -0.2 & -1.7\\
        \bottomrule
    \end{tabularx}
    \caption{Slopes of the predicted sound-localization performance paramters (PE and QE) for NH5, NH130, and NH131 from the ARI database and the average across subjects. HRTFs numerically calculated for uniform (UNI) and graded meshes (POW1 and COS2). Other details as in Tab.~\ref{tab:sphereMeshData}.}
    \label{tab:subjectsSlopes}
\end{table}


\section{Conclusions}
\label{sec:Conclusions}
A mesh preprocessing method for the numerical calculation of HRTFs, i.e., \textit{a-priori mesh grading}, was proposed. For the evaluation, the method was applied to meshes of three geometric objects with various grading functions. HRTFs were calculated based on these meshes and the results were evaluated by means of numerical errors and perceptually motivated model predictions.

When comparing the HRTFs calculated for graded meshes with those calculated for uniform meshes, the HRTF calculation for graded meshes yielded similar or even better results than for high-resolution uniform meshes in terms of numerical accuracy and in terms of the predicted sound-localization performance. Thus, graded meshes seem to be suitable for the numerical calculation of HRTFs in the full audible frequency range, further indicating that the recommendation of at least six elements per wavelength can be violated apart from the microphone area and the ipsilateral pinna.

The numerical accuracy of HRTFs calculated for various grading functions was compared including linear grading, higher-order power grading, and raised-cosine based grading. For the simple geometric object, the sphere, all grading functions showed better performance than uniform re-meshing in terms of the relative numerical error and computation time. The raised-cosine based grading and fourth-order grading functions showed the most potential. HRTFs calculated for these grading functions showed smaller errors than HRTFs calculated for the high-resolution mesh - even with less elements in the meshes. These grading functions concentrate more elements in the proximity range of the microphone area and as a consequence, less elements in total are required to obtain similar numerical errors in the numerical calculation process. Raised-cosine based grading showed the best overall performance. With approximately \num{13000} elements, the computation took only ten percent of the computation required for a high-resolution mesh, but still the relative numerical and the perceptual errors were in the range of that for the high-resolution mesh.

The proposed a-priori mesh grading algorithm was implemented as a plug-in for OpenFlipper and can be combined with the Mesh2HRTF package to calculate HRTFs. In that implementation, the proposed mesh grading algorithm reduced the calculation time and memory requirements significantly. These computational costs might be further reduced by optimizing the FMM for graded meshes and introducing a frequency-specific mesh grading in which the grading parameters are further optimized for the simulated frequency.


\section*{Acknowledgements}
\label{sec:Acknowledgements}
This study was supported by the Austrian Science Fund (FWF, projects P~24124 and I~1018~N25). The authors want to thank Andreas Reichinger (VRVis GmbH, Austria) and his team for the work on geometry acquisition required to obtain the high-resolution meshes of the HUM objects.


\bibliographystyle{elsarticle-harv}

\begin{thebibliography}{40}
\expandafter\ifx\csname natexlab\endcsname\relax\def\natexlab#1{#1}\fi
\expandafter\ifx\csname url\endcsname\relax
  \def\url#1{\texttt{#1}}\fi
\expandafter\ifx\csname urlprefix\endcsname\relax\def\urlprefix{URL }\fi

\bibitem[{Baumgartner et~al.(2013)Baumgartner, Majdak, and
  Laback}]{baumgartner_assessment_2013}
Baumgartner, R., Majdak, P., Laback, B., 2013. Assessment of {Sagittal}-{Plane}
  {Sound} {Localization} {Performance} in {Spatial}-{Audio} {Applications}. In:
  The {Technology} of {Binaural} {Listening}. Springer, Berlin, DE, pp.
  93--119.
\newline\urlprefix\url{http://dx.doi.org/10.1007/978-3-642-37762-4_4}

\bibitem[{Beranek and Mellow(2012)}]{beranek_acoustics:_2012}
Beranek, L.~L., Mellow, T.~J., Aug. 2012. Acoustics: {Sound} {Fields} and
  {Transducers}. Academic Press, Amsterdam, Netherlands.

\bibitem[{Botsch and Kobbelt(2004)}]{botsch_remeshing_2004}
Botsch, M., Kobbelt, L., 2004. A {Remeshing} {Approach} to {Multiresolution}
  {Modeling}. In: Proceedings of the 2004 {Eurographics}/{ACM} {SIGGRAPH}
  {Symposium} on {Geometry} {Processing}. {SGP} '04. ACM, New York, NY, USA,
  pp. 185--192.
\newline\urlprefix\url{http://dx.doi.org/10.1145/1057432.1057457}

\bibitem[{Burton and Miller(1971)}]{burton_application_1971}
Burton, A.~J., Miller, G.~F., Jun. 1971. The {Application} of {Integral}
  {Equation} {Methods} to the {Numerical} {Solution} of {Some} {Exterior}
  {Boundary}-{Value} {Problems}. Proceedings of the Royal Society of London. A.
  Mathematical and Physical Sciences 323~(1553), 201--210.
\newline\urlprefix\url{http://dx.doi.org/10.1098/rspa.1971.0097}

\bibitem[{Chen et~al.(2002)Chen, Chen, and Chen}]{chen_adaptive_2002}
Chen, J.~T., Chen, K.~H., Chen, C.~T., May 2002. Adaptive boundary element
  method of time-harmonic exterior acoustics in two dimensions. Computer
  Methods in Applied Mechanics and Engineering 191~(31), 3331--3345.
\newline\urlprefix\url{http://dx.doi.org/10.1016/S0045-7825(02)00214-1}

\bibitem[{Chen et~al.(2008)Chen, Waubke, and Kreuzer}]{chen_formulation_2008}
Chen, Z.-S., Waubke, H., Kreuzer, W., Jun. 2008. A formulation of the fast
  multipole boundary element method ({FMBEM}) for acoustic radiation and
  scattering from three-dimensional structures. Journal of Computational
  Acoustics 16~(2), 303--320.
\newline\urlprefix\url{http://dx.doi.org/10.1142/S0218396X08003725}

\bibitem[{Gaul et~al.(2003)Gaul, Kögl, and Wagner}]{gaul_boundary_2003}
Gaul, L., Kögl, M., Wagner, M., Feb. 2003. Boundary {Element} {Methods} for
  {Engineers} and {Scientists}: {An} {Introductory} {Course} with {Advanced}
  {Topics}. Springer, Berlin, DE.
\newline\urlprefix\url{http://dx.doi.org/10.1007/978-3-662-05136-8}

\bibitem[{Goldstein(1982)}]{goldstein_finite_1982}
Goldstein, C.~I., Feb. 1982. The finite element method with non-uniform mesh
  sizes applied to the exterior {Helmholtz} problem. Numerische Mathematik
  38~(1), 61--82.
\newline\urlprefix\url{http://dx.doi.org/10.1007/BF01395809}

\bibitem[{Gumerov et~al.(2010)Gumerov, O’Donovan, Duraiswami, and
  Zotkin}]{gumerov_computation_2010}
Gumerov, N.~A., O’Donovan, A.~E., Duraiswami, R., Zotkin, D.~N., Jan. 2010.
  Computation of the head-related transfer function via the fast multipole
  accelerated boundary element method and its spherical harmonic
  representation. The Journal of the Acoustical Society of America 127~(1),
  370--386.
\newline\urlprefix\url{http://dx.doi.org/10.1121/1.3257598}

\bibitem[{Heinrich et~al.(1996)Heinrich, Beilenhoff, Mezzanotte, and
  Roselli}]{heinrich_optimum_1996}
Heinrich, W., Beilenhoff, K., Mezzanotte, P., Roselli, L., Sep. 1996. Optimum
  mesh grading for finite-difference method. IEEE Transactions on Microwave
  Theory and Techniques 44~(9), 1569--1574.
\newline\urlprefix\url{http://dx.doi.org/10.1109/22.536606}

\bibitem[{Hunter and Pullan(2002)}]{hunter_fem/bem_2002}
Hunter, P., Pullan, A., 2002. {FEM}/{BEM} {Notes}.

\bibitem[{Jin et~al.(2014)Jin, Guillon, Epain, Zolfaghari, van Schaik, Tew,
  Hetherington, and Thorpe}]{jin_creating_2014}
Jin, C.~T., Guillon, P., Epain, N., Zolfaghari, R., van Schaik, A., Tew, A.~I.,
  Hetherington, C., Thorpe, J., Jan. 2014. Creating the {Sydney} {York}
  {Morphological} and {Acoustic} {Recordings} of {Ears} {Database}. IEEE
  Transactions on Multimedia 16~(1), 37--46.
\newline\urlprefix\url{http://dx.doi.org/10.1109/TMM.2013.2282134}

\bibitem[{Kahana and Nelson(2006)}]{kahana_numerical_2006}
Kahana, Y., Nelson, P.~A., Apr. 2006. Numerical modelling of the spatial
  acoustic response of the human pinna. Journal of Sound and Vibration
  292~(1–2), 148--178.
\newline\urlprefix\url{http://dx.doi.org/10.1016/j.jsv.2005.07.048}

\bibitem[{Kahana and Nelson(2007)}]{kahana_boundary_2007}
Kahana, Y., Nelson, P.~A., Mar. 2007. Boundary element simulations of the
  transfer function of human heads and baffled pinnae using accurate geometric
  models. Journal of Sound and Vibration 300~(3–5), 552--579.
\newline\urlprefix\url{http://dx.doi.org/10.1016/j.jsv.2006.06.079}

\bibitem[{Katz(2001{\natexlab{a}})}]{katz_boundary_2001}
Katz, B. F.~G., Nov. 2001{\natexlab{a}}. Boundary element method calculation of
  individual head-related transfer function. {I}. {Rigid} model calculation.
  The Journal of the Acoustical Society of America 110~(5), 2440--2448.
\newline\urlprefix\url{http://dx.doi.org/10.1121/1.1412440}

\bibitem[{Katz(2001{\natexlab{b}})}]{katz_boundary_2001-1}
Katz, B. F.~G., Nov. 2001{\natexlab{b}}. Boundary element method calculation of
  individual head-related transfer function. {II}. {Impedance} effects and
  comparisons to real measurements. The Journal of the Acoustical Society of
  America 110~(5), 2449--2455.
\newline\urlprefix\url{http://dx.doi.org/10.1121/1.1412441}

\bibitem[{Kreuzer et~al.(2009)Kreuzer, Majdak, and Chen}]{kreuzer_fast_2009}
Kreuzer, W., Majdak, P., Chen, Z., Sep. 2009. Fast multipole boundary element
  method to calculate head-related transfer functions for a wide frequency
  range. The Journal of the Acoustical Society of America 126~(3), 1280--1290.
\newline\urlprefix\url{http://dx.doi.org/10.1121/1.3177264}

\bibitem[{Langer et~al.(2015)Langer, Mantzaflaris, Moore, and
  Toulopoulos}]{Langer20151685}
Langer, U., Mantzaflaris, A., Moore, S.~E., Toulopoulos, I., 2015. Mesh grading
  in isogeometric analysis. Computers \& Mathematics with Applications 70~(7),
  1685 -- 1700.
\newline\urlprefix\url{http://dx.doi.org/10.1016/j.camwa.2015.03.011}

\bibitem[{Liang et~al.(1999)Liang, Chen, and Yang}]{liang_error_1999}
Liang, M.~T., Chen, J.~T., Yang, S.~S., Mar. 1999. Error estimation for
  boundary element method. Engineering Analysis with Boundary Elements 23~(3),
  257--265.
\newline\urlprefix\url{http://dx.doi.org/10.1016/S0955-7997(98)00086-1}

\bibitem[{Luebke(2001)}]{luebke_developers_2001}
Luebke, D., May 2001. A developer's survey of polygonal simplification
  algorithms. IEEE Computer Graphics and Applications 21~(3), 24--35.
\newline\urlprefix\url{http://dx.doi.org/10.1109/38.920624}

\bibitem[{Macpherson and Middlebrooks(2002)}]{macpherson_listener_2002}
Macpherson, E.~A., Middlebrooks, J.~C., May 2002. Listener weighting of cues
  for lateral angle: {The} duplex theory of sound localization revisited. The
  Journal of the Acoustical Society of America 111~(5), 2219--2236.
\newline\urlprefix\url{http://dx.doi.org/10.1121/1.1471898}

\bibitem[{Majdak et~al.(2007)Majdak, Balazs, and Laback}]{majdak_multiple_2007}
Majdak, P., Balazs, P., Laback, B., Jul. 2007. Multiple {Exponential} {Sweep}
  {Method} for {Fast} {Measurement} of {Head}-{Related} {Transfer} {Functions}.
  Journal of the Audio Engineering Society 55~(7/8), 623--637.
\newline\urlprefix\url{http://www.aes.org/e-lib/browse.cfm?elib=14190}

\bibitem[{Majdak et~al.(2010)Majdak, Goupell, and Laback}]{majdak_3-d_2010}
Majdak, P., Goupell, M.~J., Laback, B., Feb. 2010. 3-{D} {Localization} of
  {Virtual} {Sound} {Sources}: {Effects} of {Visual} {Environment}, {Pointing}
  {Method}, and {Training}. Attention, perception \& psychophysics 72~(2),
  454--469.
\newline\urlprefix\url{http://dx.doi.org/10.3758/APP.72.2.454}

\bibitem[{Marburg(2002)}]{marburg_six_2002}
Marburg, S., Mar. 2002. Six boundary elements per wavelength: {Is} that enough?
  Journal of Computational Acoustics 10~(1), 25--51.
\newline\urlprefix\url{http://dx.doi.org/10.1142/S0218396X02001401}

\bibitem[{Mehrgardt and Mellert(1977)}]{mehrgardt_transformation_1977}
Mehrgardt, S., Mellert, V., Jun. 1977. Transformation characteristics of the
  external human ear. The Journal of the Acoustical Society of America 61~(6),
  1567--1576.
\newline\urlprefix\url{http://dx.doi.org/10.1121/1.381470}

\bibitem[{Middlebrooks(1999)}]{middlebrooks_individual_1999}
Middlebrooks, J.~C., Sep. 1999. Individual differences in external-ear transfer
  functions reduced by scaling in frequency. The Journal of the Acoustical
  Society of America 106~(3), 1480--1492.
\newline\urlprefix\url{http://dx.doi.org/10.1121/1.427176}

\bibitem[{Morimoto and Aokata(1984)}]{morimoto_localization_1984}
Morimoto, M., Aokata, H., 1984. Localization cues of sound sources in the upper
  hemisphere. The Journal of the Acoustical Society of Japan~(5), 166--173.
\newline\urlprefix\url{http://dx.doi.org/10.1250/ast.5.165}

\bibitem[{Morse and Ingard(1986)}]{morse_theoretical_1986}
Morse, P.~M., Ingard, K.~U., 1986. Theoretical acoustics. Princeton University
  Press, Princeton, NJ.

\bibitem[{Möbius and Kobbelt(2012)}]{mobius_openflipper:_2012}
Möbius, J., Kobbelt, L., 2012. {OpenFlipper}: {An} {Open} {Source} {Geometry}
  {Processing} and {Rendering} {Framework}. In: Proceedings of the 7th
  {International} {Conference} on {Curves} and {Surfaces}. Springer-Verlag,
  Berlin, DE, pp. 488--500.
\newline\urlprefix\url{http://dx.doi.org/10.1007/978-3-642-27413-8_31}
\vfill\break
\bibitem[{Møller(1992)}]{moller_fundamentals_1992}
Møller, H., 1992. Fundamentals of binaural technology. Applied Acoustics
  36~(3–4), 171--218.
\newline\urlprefix\url{http://dx.doi.org/10.1016/0003-682X(92)90046-U}

\bibitem[{Møller et~al.(1995)Møller, Sørensen, Hammershøi, and
  Jensen}]{moller_head-related_1995}
Møller, H., Sørensen, M.~F., Hammershøi, D., Jensen, C.~B., May 1995.
  Head-{Related} {Transfer} {Functions} of {Human} {Subjects}. Journal of the
  Audio Engineering Society 43~(5), 300--321.
\newline\urlprefix\url{http://www.aes.org/e-lib/browse.cfm?elib=7949}

\bibitem[{Reichinger et~al.(2013)Reichinger, Majdak, Sablatnig, and
  Maierhofer}]{reichinger_evaluation_2013}
Reichinger, A., Majdak, P., Sablatnig, R., Maierhofer, S., Jun. 2013.
  Evaluation of {Methods} for {Optical} 3-{D} {Scanning} of {Human} {Pinnas}.
  In: 2013 {International} {Conference} on 3D {Vision} - 3DV 2013. pp.
  390--397.
\newline\urlprefix\url{http://dx.doi.org/10.1109/3DV.2013.58}

\bibitem[{Rui et~al.(2013)Rui, Yu, Xie, and Liu}]{rui_calculation_2013}
Rui, Y., Yu, G., Xie, B., Liu, Y., May 2013. Calculation of {Individualized}
  {Near}-{Field} {Head}-{Related} {Transfer} {Function} {Database} {Using}
  {Boundary} {Element} {Method}. In: Proceedings of the 134th {Convention} of
  the Audio Engineering Society. Rome, IT.

\bibitem[{Søndergaard and Majdak(2013)}]{sondergaardmajdak2013}
Søndergaard, P., Majdak, P., 2013. The {Auditory} {Modeling} {Toolbox}. In:
  Blauert, J. (Ed.), The {Technology} of {Binaural} {Listening}. Springer,
  Berlin, DE, pp. 33--56.
\newline\urlprefix\url{http://dx.doi.org/10.1007/978-3-642-37762-4_2}

\bibitem[{Treeby and Pan(2009)}]{treeby_practical_2009}
Treeby, B.~E., Pan, J., Nov. 2009. A practical examination of the errors
  arising in the direct collocation boundary element method for acoustic
  scattering. Engineering Analysis with Boundary Elements 33~(11), 1302--1315.
\newline\urlprefix\url{http://dx.doi.org/10.1016/j.enganabound.2009.06.005}

\bibitem[{Walsh and Demkowicz(2003)}]{walsh_hp_2003}
Walsh, T., Demkowicz, L., Jan. 2003. hp {Boundary} element modeling of the
  external human auditory system––goal-oriented adaptivity with multiple
  load vectors. Computer Methods in Applied Mechanics and Engineering
  192~(1–2), 125--146.
\newline\urlprefix\url{http://dx.doi.org/10.1016/S0045-7825(02)00536-4}

\bibitem[{Wightman and Kistler(1989)}]{wightman_headphone_1989}
Wightman, F.~L., Kistler, D.~J., Feb. 1989. Headphone simulation of free-field
  listening. {I}: Stimulus synthesis. The Journal of the Acoustical Society of
  America 85~(2), 858--867.
\newline\urlprefix\url{http://dx.doi.org/10.1121/1.397557}

\bibitem[{Ziegelwanger and Majdak(2014)}]{ziegelwanger_modeling_2014}
Ziegelwanger, H., Majdak, P., Mar. 2014. Modeling the direction-continuous
  time-of-arrival in head-related transfer functions. The Journal of the
  Acoustical Society of America 135~(3), 1278--1293.
\newline\urlprefix\url{http://dx.doi.org/10.1121/1.4863196}

\bibitem[{Ziegelwanger et~al.(2015{\natexlab{a}})Ziegelwanger, Majdak, and
  Kreuzer}]{ziegelwanger_mesh2hrtf:_2015}
Ziegelwanger, H., Majdak, P., Kreuzer, W., 2015{\natexlab{a}}. Mesh2hrtf:
  {Open}-source software package for the numerical calculation of head-related
  transfer functions. In: Proceedings of the 22nd {International} {Congress} of
  {Sound} and {Vibration}. Florence, IT.

\bibitem[{Ziegelwanger et~al.(2015{\natexlab{b}})Ziegelwanger, Majdak, and
  Kreuzer}]{ziegelwanger_numerical_2015}
Ziegelwanger, H., Majdak, P., Kreuzer, W., Jul. 2015{\natexlab{b}}. Numerical
  calculation of listener-specific head-related transfer functions and sound
  localization: {Microphone} model and mesh discretization. The Journal of the
  Acoustical Society of America 138~(1), 208--222.
\newline\urlprefix\url{http://dx.doi.org/10.1121/1.4922518}

\end{thebibliography}


\end{document}